%
\magnification=\magstep1
\input amstex
\UseAMSsymbols
\input pictex
\NoBlackBoxes
\font\kl=cmr8  \font\gross=cmbx10 scaled\magstep1 \font\abs=cmcsc10
   \font\rmk=cmr8    \font\itk=cmti8  \font\ttk=cmtt8

   \newcount\notenumber
   
   \def\note{\advance\notenumber by 1 
       \plainfootnote{$^{\the\notenumber}$}}  

\def\mod{\operatorname{mod}}

\def\Mat{\operatorname{Mat}}
\def\Hom{\operatorname{Hom}}
\def\ind{\operatorname{ind}}
\def\Tor{\operatorname{Tor}}
\def\End{\operatorname{End}}
\def\Ext{\operatorname{Ext}}
\def\rad{\operatorname{rad}}

\def\add{\operatorname{add}}

\def\abs{\operatorname{abs}}
\def\supp{\operatorname{supp}}

\def\bdim{\operatorname{\bold{dim}}}

\def\arr#1#2{\arrow <1.5mm> [0.25,0.75] from #1 to #2}

   


\centerline{\gross Cluster-concealed algebras}     
		   \bigskip
\centerline{Claus Michael Ringel}
		  \bigskip\medskip
\plainfootnote{}
{\rmk 2000 \itk Mathematics Subject Classification. \rmk 
Primary 
        16D90, 
        16G60. 
Secondary:
        16G20, 
        16G70. 
}

	\bigskip\bigskip
{\narrower\narrower\narrower \kl Abstract. The cluster-tilted algebras have been introduced
by Buan, Marsh and Reiten, they are the endomorphism rings of cluster-tilting objects $\ssize T$ 
in cluster categories; we call such an algebra cluster-concealed in case $\ssize T$ is
obtained from a preprojective tilting module. For example, all representation-finite
cluster-tilted algebras are cluster-concealed. If $\ssize C$ is 
a representation-finite cluster-tilted algebra, then the indecomposable $\ssize C$-modules
are shown to be determined by their dimension vectors.
For a general cluster-tilted algebra $\ssize C$, we are going to describe the dimension vectors
of the indecomposable $\ssize C$-modules in terms of the root system of a
quadratic form. The roots may have both positive and
negative coordinates and we have to take absolute values. 
\par}
	\bigskip\bigskip

Let $k$ be an algebraically closed field. For any finite-dimensional $k$-algebra $R$, 
we consider its Grothendieck group 
$K_0(R)$ of (finitely generated) $R$-modules modulo exact sequences: it is the
free abelian group with basis the set of isomorphism classes of simple $R$-modules.
Using this basis, we identify $K_0(R)$ with $\Bbb Z^n$,
where $n$ is the number of isomorphism classes of simple $R$-modules.
For any $R$-module $N$, we denote by $\bdim N$ the corresponding element in $K_0(R)$.
With respect to our identification $K_0(R) = \Bbb Z^n$, the coefficients of $\bdim N$ 
are just the Jordan-H\"older multiplicities of $N$ and the set of simple $R$-modules 
which occur as composition factors of $N$ will be called the {\it support} of $N$
and denotes by $\supp N$.

    \medskip
Throughout the paper, $A$ will denote a finite-dimensional hereditary $k$-algebra.
Recall that a $k$-algebra $B$ is said to be {\it tilted} 
provided $B$ is the endomorphism ring
of a tilting $A$-module $T$, where $A$ is a finite-dimensional hereditary
algebra.  If $B$ is a tilted
algebra, one may consider the corresponding trivial extension
algebra $B^c = B \ltimes I$, where $I$ is the $B$-$B$-bimodule
$I = \Ext^2_B(DB,B)$, with $D = \Hom(-,k)$  the $k$-duality. The algebras
of the form $B^c$ are called the {\it cluster tilted} algebras; 
this is not the
original definition as given by
Buan, Marsh and Reiten [BMR], but it is an equivalent
one, due to Zhu [Z] and Assem, Br\"ustle and Schiffler [ABS].
The definition shows that 
$B$ is both a subalgebra as well as a factor algebra of $B^c$,
and that the $C$-modules $N$ with $IN = 0$ are just the $B$-modules.
    \medskip
    
{\bf Theorem 1.} {\it Let $C$ be a representation-finite cluster-tilted algebra.
If $N,N'$ are indecomposable $C$-modules with $\bdim N = \bdim N'$, then $N$ and
$N'$ are isomorphic.}
     \medskip
After finishing a first version of this paper, we were informed about a parallel 
investigation by Geng and Peng [GP] which gives a different proof of this result 
using mutations.  The paper by Geng and Peng also outlines the link to 
cluster algebras and shows in which way Theorem 1 settles
a conjecture of Fomin and Zelevinsky concerning cluster variables.
  \bigskip
The next result will provide a description of the set of dimension vectors $\bdim N$ in
$K_0(C)$ with $N$ indecomposable. This will be done in a more general setting. 
Recall  that a tilted algebra $B$ is said to be {\it concealed} 
provided $B$ is the endomorphism ring
of a preprojective tilting $A$-module.
If $B$ is a concealed algebra, then we will say that $B^c$ is a 
{\it cluster-concealed}
algebra. Of course, representation-finite cluster-tilted algebras are cluster-concealed
algebras, but there are also many cluster-concealed algebras which are tame or wild.
	  \medskip
A famous theorem of Kac asserts that the dimension vectors of the indecomposable
$A$-modules are just the positive roots of the (Kac-Moody) root system $\Phi_A$
in  $K_0(A)$ corresponding to the underlying graph of the quiver of $A$.
Note that $q_A(x) \le 1$ for any $x\in \Phi_A$, 
here $q_A$ is the Euler form on $K_0(A)$ (the definition
will be recalled in section 10).
     \medskip
The reason for us to
exhibit cluster-tilted algebras as $B^c = B\ltimes I$ is that this allows
to identify the Grothendieck groups $K_0(B^c)$ and $K_0(B)$.
Let $T$ be a tilting $A$-module with endomorphism ring $B$, and let $q_B$ be the 
Euler form of $B$ on $K_0(B)$. Since we identify 
the Grothendieck groups $K_0(B) = K_0(B^c)$, we can apply $q_B$ to 
the dimension vectors of the indecomposable $B^c$-modules; this is what we
will do. On the other hand,
consider the tilting functor $G = \Hom_A(T,-)\:\mod A \to \mod B$.
Let $T_1,\dots,T_n$ be indecomposable direct summands of $T$, one from each
isomorphism class. Then $\bdim  T_1,\dots, \bdim  T_n$ is a basis of $K_0(A)$,
whereas $\bdim  G(T_1),\dots,\bdim  G(T_n)$ is a basis of $K_0(B),$ and 
we denote by $g\:K_0(A) \to K_0(B)$ the linear bijection such that
$g(\bdim  T_i) = \bdim  G(T_i)$, for $1 \le i \le n$. We set $\Phi_B = g(\Phi_A).$
If $x\in \Phi_A$, then it is well-known that $x$ or $-x$ belongs to $\Bbb N$, but,
usually, $\Phi_B$ will contain elements for which some coefficients are
positive, and some negative. For any element $x = (x_1,\dots,x_n)\in \Bbb Z^n$, 
let $\abs x = (|x_1|,\dots,|x_n|).$ We use this function $\abs$ in order to 
attach to any vector $x\in \Phi_B$ an element in $\Bbb N^n$. 
       \bigskip
{\bf Theorem 2.} {\it Let $B$ be a concealed algebra and $C = B^c$ the
corresponding cluster-concealed algebra,

\item{\rm(a)} The dimension vectors of the indecomposable $C$-modules
are precisely the vectors $\abs x$ with $x\in \Phi_B.$ 
\item{\rm(b)} If $Z$ is an indecomposable
$C$-module, then $q_B(\bdim N) \le 1$ if and only if $Z$ is a $B$-module; 
if $N$ is not a $B$-module, then  
$q_B(\bdim N)$ is an odd integer (greater than $2$).
\item{\rm(c)} If $N$ is an indecomposable
$C$-module which is not a $B$-module, then $\End(N) = k.$} 
	   \bigskip
	   
The 
special case when $A$ is of type $\Bbb A_n$ has been considered already in the thesis
of Parsons [P], using a different approach. 
As Robert Marsh has pointed out, some further cases have been considered by
Parsons and him in this way but this work is still ongoing.
	   \bigskip
\vfill\eject
{\bf Remarks.}
     \medskip
{\bf  1. The quadratic form $q_B$.} 
We want to stress that the quadratic form $q_B$ used here in order to
deal with $B^c$-modules depends on the choice of $B$: it is the Euler form
for $B$, not for $B^c$ (there may not even exist a Euler form for $B^c$,
since usually the global dimension of $B^c$ is infinite). Also, for a given 
cluster-concealed algebra $C$, there usually will exist several concealed algebras $B$
with $C = B^c$ and then we will obtain different quadratic forms $q_B$ on
$K_0(C),$ see Example 13.2. 
	  \medskip
	  
{\bf 2. About the proof.} 
The upshot of our investigation is Proposition 4. 
A direct proof of this result (as well as 
a generalization to tilting modules which are neither preprojective nor
preinjective) would be of interest. Our proof
invokes a second quadratic form $r_E$ which 
concerns the extension behavior of two torsion pairs. Here, we
deal with what Drozd [D] has called $E$-matrices where $E$ is a bimodule. The 
preprojectivity of $T$ is used in order to show that the corresponding 
categories of $E$-matrices are essentially directed: 
according to de la Pe\~na and Simson [DP] this then implies that the
indecomposable objects correspond bijectively to the positive roots
of the corresponding quadratic forms. But actually, these quadratic forms 
coincide, in this way we obtain the required bijection. We should stress
that the equality of the quadratic forms used follows from the fact that
the bimodules $E$ which arise can be written as
$\Ext^1(\Cal F,\Cal G)$ and as $\Hom(\Cal G,\tau \Cal F)$, respectively, and, 
of course, one of the
basic results of Auslander-Reiten theory asserts that these bimodules are dual to each
other.
	
A second ingredient of our proof is the following separation property (see section 2): 
Let $T$ be a preprojective tilting module, and $M$ an indecomposable
$A$-module. Then we show that the $B$-modules
$G(M) = \Hom_A(T,M)$ and  $F(M) = \Ext^1(T,M)$
have disjoint supports. This property is the reason for the appearance
of absolute values in Theorem 2. In case $T$ has an indecomposable regular
direct summand, the separation 
property no longer holds, see Proposition 7. Thus, one cannot expect that a
generalization of the main theorem for arbitrary cluster-tilted algebras 
will use the vectors $\abs x$ with $x\in \Phi_B.$ 
     
Invariants such as quadratic forms or root systems have often been
used in order
to obtain a classification of the indecomposable $R$-modules, for an algebra $R$.
Usually, one starts to guess all these objects, then one shows that they are 
pairwise non-isomorphic and that all 
the indecomposable $R$-modules have been obtained in this way, 
and finally, one tries to describe the structure of the module category.
In our case of dealing with a cluster-tilted algebra, 
the procedure is completely reversed:
the module category is known from the beginning, but one is lacking sufficient 
information concerning the modules themselves.
	    \medskip
	    
{\bf 3. The relevance of cluster-concealed algebras.} 
The importance of the concealed algebras should be mentioned here. 
The tame ones have been classified by Happel and Vossieck [HV], and are used 
by the Bongartz criterion for determining whether a finite-dimensional algebra
is representation-finite or not. The list of all
the possible ``frames'' of tame concealed algebras can be found in several books 
and papers, the corresponding cluster algebras have been considered by Seven [S],
the relationship has been discussed in [BRS].
Wild concealed algebras have been considered by Unger [U].

     \bigskip
{\bf Acknowledgment.} 
The results presented in the paper have been obtained during a stay at the
Newton Institute in Cambridge, February and March 2009. The author is grateful
to the Institute for its hospitality. Also, the author is
indebted to Idun Reiten and the referee for a careful reading of a first version. 

	 \bigskip\bigskip
{\bf 1. Notation.} 
     \medskip

If $R$ is a finite-dimensional $k$-algebra, the modules considered usually
will be finite-dimensional left modules, and
homomorphisms will be written at the opposite side of the scalars. 
Let $\mod R$ be
the category of $R$-modules, and $\ind R$ a set 
of indecomposable $R$-modules, one from each isomorphism class, or also the
corresponding full subcategory of $\mod R$.
	      \medskip
We denote by $A$ a finite-dimensional hereditary $k$-algebra which
is connected. Let $T$ be a tilting $A$-module with endomorphism ring $B$.
As usually, we let
$$
	\align
\Cal F = \Cal F(T)  &= \{ M \in \ind A \mid \Hom(T,M) = 0\}, \cr
\Cal G = \Cal G(T)  &= \{ M \in \ind A \mid \Ext^1(T,M) = 0\}, \cr
\Cal X = \Cal X(T)  &= \{ M \in \ind B \mid T\otimes_B M = 0\}, \cr
\Cal Y = \Cal Y(T)  &= \{ M \in \ind B \mid \Tor_1(T,M) = 0\}.
     \endalign
$$
The pair $(\Cal F,\Cal G)$ is a torsion pair in $\mod A$. Given an $A$-module $M$,
we denote its torsion submodule by $tM$. The pair 
$(\Cal Y,\Cal X)$ is a torsion pair in $\mod B$ which is even split.
       \medskip
Tilting theory asserts that the functor $G = \Hom(T,-)$ gives an equivalence
$\Cal G(T) \to \Cal Y(T)$ and that the functor $F = \Ext^1(T,-)$ gives an equivalence
$\Cal F(T) \to \Cal X(T)$.
We should stress that for any $A$-module $M$
$$
\align
 G(M) &= \Hom(T,M) = \Hom(T,tM) = G(tM),\cr
 F(M) &= \Ext^1(T,M) = \Ext^1(T,M/tM) = F(M/tM).
\endalign
$$
	\medskip
In this paper, the main interest will concern
$$
 \Cal M = \Cal M(T)  = \{ M \in \ind A \mid \Hom(T,M) \neq 0,\ \Ext^1(T,M) \neq 0\},
$$
as well as 
$$
 \Cal N = \Cal N(B) = \{ N \in \ind B^c \mid IN \neq 0\},
$$
these are the indecomposable $B^c$-modules which are not $B$-modules.
      \bigskip\bigskip

{\bf 2. $\Cal M(T)$ for $T$ preprojective.}
     \medskip
Recall the following: Given a $B$-module $N$, its {\it support} is the 
set of isomorphism classes of the
simple $B$-modules which occur as  composition factors of $N$. The indecomposable
projective $B$-modules are of the form $G(T_i)$, thus the simple $B$-modules are
indexed by the natural numbers $1,2,\dots,n$. Tilting theory shows that for any
$A$-module $M$, the support $\supp G(M)$ of $G(M)$ is the set of indices $i$
such that $\Hom(T_i,M) \neq 0$, and that the support $\supp F(M)$ of $F(M)$ is the set of
indices $i$ such that $\Ext^1(T_i,M) \neq 0$.
	\medskip
{\bf Lemma 1.} {\it Let $M, M'$ be indecomposable $A$-modules. Assume that there is
an index $i$ in the intersection of $\supp F(M')$ and $\supp G(M)$ such that $T_i$ is
preprojective or preinjective.  Then $\supp F(M)$ and $\supp G(M')$
do not intersect.}
   \medskip
Proof: Since $i\in\supp F(M')\cap \supp G(M)$, we have
$$
 \Ext^1(T_i,M') \neq 0 \quad\text{and}\quad \Hom(T_i,M) \neq 0.
$$ 
Assume that there is  $j$ in the intersection of  $\supp F(M)$ and $\supp G(M')$, thus
$$
\Ext^1(T_j,M) \neq 0 \quad\text{and}\quad \Hom(T_j,M') \neq 0.
$$ 
Note that $\Ext^1(T_i,M') = D\Hom(M',\tau T_i)$ and 
$\Ext^1(T_j,M) = D\Hom(M,\tau T_j)$. Thus, we obtain a proper cycle 
$$
 T_i \preccurlyeq M \preccurlyeq \tau T_j \prec T_j \preccurlyeq M' \preccurlyeq \tau T_i 
 \prec T_i
$$
whereas $T_i$ is preprojective or preinjective, thus directing. 
	\medskip
The case  $M = M'$ yields the following corollary:
    \medskip
{\bf Corollary.} {\it Let $M$ be an indecomposable $A$-module. Assume that there is
an index $i$ in the intersection of $\supp F(M')$ and $\supp G(M)$, then $T_i$ is neither
preprojective nor preinjective.}
	      \medskip
Thus we have the following
     \medskip
{\bf Separation Property.} {\it If $T$ is a preprojective 
tilting module, and $M$
is indecomposable, then the supports of $F(M)$ and $G(M)$ are
disjoint.}

	\bigskip
For the remainder of this section, we assume that $T$ is a preprojective tilting module.
Let $\Cal D$ be the set
of predecessors of the modules $\tau T_i$, where $T_i$ is an
indecomposable  direct summand of $T$.
		\medskip
{\bf Lemma 2.} {\it We have $\Cal F \subseteq \Cal D$.}
     \medskip
Proof: Let $X$ be in $\Cal F$. Since
$F\:\Cal F \to \Cal X$ is an equivalence, we have $\Ext^1(T,X) = FX \neq 0$,
thus $0 \neq \Ext^1(T_i,X) = D\Hom(X,\tau T_i)$ for some $i$. This means that
$X$ is a predecessor of $\tau T_i$, thus is in $\Cal D$. 
    \bigskip
{\bf Lemma 3.} {\it If $M$ belongs to $\Cal M$, then
any indecomposable direct summand of $M,$ $tM,$ $M/tM$ belongs to $\Cal D.$}
    \medskip
Proof: For any module $M$, the indecomposable direct summands of $M/tM$
belong to $\Cal F$, thus to $\Cal D$, according to Lemma 2. 
For
$M\in \Cal M$, the factor module $M/tM$ is non-zero, thus $M$ has an
indecomposable factor module $M'$ which belongs to $\Cal F$, thus
also $M$ belongs to $\Cal D$. But then also any indecomposable summand of any
submodule of $M$ belongs to $\Cal D$. In particular, any indecomposable summand
of $tM$ belongs to $\Cal D.$
   \medskip

{\bf Lemma 4.} {\it 
If $M \in \Cal M,$ then $\Hom(M/tM,tM) = 0$ and
both $tM$ and $M/tM$ have no self-extensions.}
     \medskip
Proof: Since $M$ is indecomposable and preprojective, its endomorphism ring is $k$.
Any non-zero homomorphism $M/tM \to tM$ would yield a non-zero nilpotent endomorphism
of $M$, which is impossible. 

Let $M'$ be an indecomposable direct summand of $tM$.
Then we get the following exact sequence
$$
  \Hom(M',M/tM) \to \Ext^1(M',tM) \to \Ext^1(M',M).
$$
Here, the first term is zero, since $M'$ is torsion, and $M/tM$ torsion-free. Also, the 
last term is zero, since $M'$ is a predecessor of $M$, thus there cannot be
a proper path from $M$ to $M'$. Thus $\Ext^1(M',tM) = 0$ and therefore $\Ext^1(tM,tM) = 0.$

Similarly, let $M'$ be an indecomposable direct summand of $M/tM$. There is the
following exact sequence
$$
 \Hom(tM,M') \to \Ext^1(M/tM,M') \to \Ext^1(M,M').
$$
Again, the first term is zero, since $tM$ is torsion and $M'$ is torsion-free.
The last term is zero, since $M'$ is a successor of $M$, thus there cannot be a
proper path from $M'$ to $M$. Therefore, $\Ext^1(M/tM,M') = 0,$ thus
$\Ext^1(M/tM,M/tM) = 0.$
		   
		   \bigskip\bigskip

{\bf 3. The matrix category of a bimodule.}
     \medskip
Given two additive categories $\Cal A$ and $\Cal B$, an $\Cal A$-$\Cal B$-bimodule
${}_{\Cal A}E_{\Cal B}$ is by definition 
a bilinear functor $\Cal A^{\text{op}}\times \Cal B \to\mod k.$ 
Given such an $\Cal A$-$\Cal B$-bimodule $E = {}_{\Cal A}E_{\Cal B}$,
let $\Mat E$ be the category of {\it $E$-matrices}
as introduced by Drozd [D]: its objects are triples $(A,B,m)$, where $A$ is an object of
$\Cal A$, $B$ is an object of $\Cal B$ and $m\in E(A,B)$. A morphism
$(A,B,m) \to (A',B',m')$ is a pair $(\alpha,\beta)$, where $\alpha\:A \to A'$
and $\beta\:B \to B'$ are morphisms in $\Cal A$, and  $\Cal B$ respectively,
such that $m\beta = \alpha m'.$
     \bigskip
Given a bimodule ${}_{\Cal A}E_{\Cal B}$, one may introduce a
quadratic form $r_E$ as follows: it is defined on the direct sum of the
Grothendieck groups $K(\Cal A,\oplus)$ and $K(\Cal B,\oplus)$. If $X$ is an
object in $\Cal A$, and $Y$ an object in $\Cal B$, then one calls the pair
$(X,Y)\in K(\Cal A,\oplus) \oplus K(\Cal B,\oplus)$ a {\it coordinate vector} and
one puts
$$
 r_E((X,Y)) = \dim\End_{\Cal A}(X) + \dim \End_{\Cal B}(Y) - \dim E(X,Y),
$$
and this extends in a unique way to a quadratic form on 
$K(\Cal A,\oplus) \oplus K(\Cal B,\oplus)$. This quadratic form can be
presented by drawing a graph with two kinds of edges, say solid ones and
dotted ones. 

We recall the following: Assume there is given 
a free abelian group $K$ with a fixed basis
$\Cal B$ and a quadratic form $q$ on $K$ with integer values
such that $q(b) = 1$ for all $b\in \Cal B$, then we draw the following graph:
its vertices are the elements of $\Cal B$; and for $b\neq b'$ in $\Cal B$, we consider
$c_{bb'} = q(b+b')$. If $c_{bb'}$ is negative, then 
we draw $-c_{bb'}$ solid edges between $b$ and $b'$, otherwise we draw $c_{bb'}$
dotted edges between $b$ and $b'$. In the case of the quadratic form $r_E$ on 
$K = K(\Cal A,\oplus) \oplus K(\Cal B,\oplus)$, we take as $\Cal B$ the 
set of indecomposble objects in $\Cal A$ and $\Cal B$ and observe that 
the required condition $r_E(b) = 1$ for $b\in \Cal B$ is satisfied.
Thus, in our case the vertices are the isomorphism classes of the indecomposable
objects in $\Cal A$ and $\Cal B$, there are solid edges between vertices for
$\Cal A$ and $\Cal B$, and there are dotted edges between the vertices for $\Cal A$
as well as between the vertices for $\Cal B$; in this way, the graph is {\it bipartite}.
The number of edges is as follows: for indecomposable objects $X,X'$ both in $\Cal A$, 
or both in $\Cal B$, the number of dotted edges is $\dim\rad(X,X') + \dim \rad(X',X)$,
whereas the number of solid edges between an isomorphism class $A$ in $\Cal A$ and an
isomorphism class $B$ in $\Cal B$ is $\dim E(A,B)$. 
	    \bigskip
A Krull-Remak-Schmidt category with finitely many isomorphism classes of indecomposable
objects is said to be {\it directed}, provided the endomorphism rings of 
all the indecomposable objects are division rings and there is a
total ordering $\prec$ on the set of isomorphism classes of
indecomposable objects such that $\Hom(M,M') = 0$ for $M' \prec M.$
	       \medskip
There is the following important result of de la Pe\~na and Simson ([DS],
Propositions 1.1, 4.2 and 4.13):
{\it If the category $\Mat E$ is directed, then $r_E$ is weakly positive, and the use
of coordinate vectors provides a bijection 
between the indecomposable objects in $\Mat E$
and the positive roots of $r_E$.}
    \medskip
When dealing with a bimodule $E = {}_{\Cal A}E_{\Cal B}$, we sometimes will
write $E = (\Cal A,E,\Cal B)$ in order to stress the categories $\Cal A, \Cal B$.
As usually, we may factor our the annihilators, thus we may consider
$\overline E = (\overline{\Cal A},\overline E,\overline{\Cal B})$, 
where $\overline{\Cal A}$
is obtained from $\Cal A$ by factoring out the ideal of all maps $\alpha\in \Cal A$
such that $E(\alpha,1_B) = 0$ for all objects $B$ in $\Cal B$, where similarly
$\overline{\Cal B}$
is obtained from $\Cal B$ by factoring out the ideal of all maps $\beta\in \Cal B$
such that $E(1_B,\beta) = 0$ for all objects $A$ in $\Cal A$, and 
$\overline E(\overline\alpha,\overline\beta) = E(\alpha,\beta).$ We say that
$E = (\Cal A,E,\Cal B)$ is {\it essentially directed}, provided that 
$\overline E = (\overline{\Cal A},\overline E,\overline{\Cal B})$ is directed.
	   \bigskip\bigskip
{\bf 4. The category $\Mat \Ext^1(\Cal F,\Cal G)$ for $T$ preprojective.}
     \medskip
The first bimodule to be considered is $\Ext^1(\Cal F,\Cal G)$: here, $\Cal A = \Cal F$,
$\Cal B = \Cal G$ and the functor is $E = \Ext^1(-,-).$ 
      \medskip
{\bf Proposition 1.} {\it If $M$ is an $A$-module, let $\eta(M) = (M/tM,tM;\epsilon)$, 
where $\epsilon$ is the
equivalence class of the canonical exact sequence $0 \to tM \to M \to M/tM \to 0$. This
defines a functor $\eta\: \mod A \to \Mat\Ext^1(\Cal F,\Cal G)$ which 
is full and dense and its kernel is the ideal generated by all maps $\Cal F \to \Cal G$.}
   \medskip 
Proof: well-known and easy (see the appendix).   
       \bigskip

{\bf Corollary.} {\it Assume that $T$ is preprojective.
The category $\Mat \Ext^1(\Cal F,\Cal G)$ is essentially
directed.} 
	   \medskip
Proof: The functor $\eta$ maps $\Cal D$ onto $\Mat \Ext^1(\Cal F,\Cal G\cap \Cal D)$.
Now, $\Cal D$ is directed, thus also $\Mat \Ext^1(\Cal F,\Cal G\cap \Cal D)$ is directed. 
     \bigskip\bigskip
     
{\bf 5. The bimodule $I$ and the algebra $B_2$.}
     \medskip
We consider the $B$-$B$-bimodule $I = \Ext_B^2(DB,B)$. According to
[ABS] (see also [R]), this bimodule can be identified with 
$\Ext^1_A(T,\tau^{-1}T) = F(\tau^{-1}T) $.
 
 \bigskip
{\bf Lemma 5.} {\it For $X \in \Cal F$, the $B$-modules 
$\Hom_B(I,F(X))$  and $\Hom_A(T,\tau X)$ are isomorphic.}
		  \medskip
Proof: 
$$
\align
 \Hom_A(T,\tau X)  &\cong \Hom_A(\tau^{-1}T,X) \cr
 &\cong \Hom_B(F(\tau^{-1}T),F(X)) = \Hom_B(I,F(X)).
\endalign
$$
First, we have used that $\tau^{-1}$ and $\tau$ are adjoint, and then the fact that
$F$ yields a bijection $\Hom_A(M,X) \to \Hom_B(F(M),F(X)),$ for any $A$-module $M$,
since $X$ belongs to $\Cal F$ (note that this bijection for $M = \tau^{-1}T$ 
is $B$-linear, since 
$\tau^{-1}T$ is an $A$-$B$-bimodule).

	     \bigskip
{\bf Lemma 6.} {\it 
For any $B$-module $N$, the module $I\otimes_B N$ belongs to
$\add \Cal X,$ and the module $\Hom_B(I,N)$ to $\add \Cal Y$.}
      \medskip
Proof: Note that $I = \Ext^1(T,\tau^{-1}T)$ belongs to $\add \Cal X$.
Let $p\:B^t \to N$ be a free cover of $N$, then $I\otimes p\: I^t = I\otimes_B B^t \to
I \otimes_B N$ is surjective. Thus, with $I$ also $I\otimes_B N$ belongs to 
$\add\Cal X$.

In order to show that $\Hom_B(I,N)$ belongs to $\add \Cal Y$, decompose
 $N = N'\oplus N''$ with $N'\in \add \Cal X$ and $N'' \in \add\Cal Y.$ Then 
$\Hom_B(I,N) = \Hom_B(I,N'),$ thus we can assume that $N\in \add\Cal X$.
The previous lemma asserts that $\Hom_B(I,N)$ is isomorphic (as a $B$-module)
to $\Hom_A(T,\tau M)$,  with $M\in \add\Cal F$. But $\Hom_A(T,\tau M) = 
\Hom_A(T,t(\tau M))$ belongs to $\add \Cal Y$.
		\bigskip
{\bf Remark.} Lemma 6 implies the (well-known) fact 
that $I\otimes_B I = 0.$ Namely, consider
the adjoint map $\delta\:I \to \Hom_B(I,I\times I)$ of the identity map. 
We know that $I\in \add \Cal X$, whereas $\Hom_B(I,I\times I)$ belongs to $\add
\Cal Y$. Thus $\delta = 0$, and therefore $I\otimes I = 0.$ --- 
Note that $I\otimes I = 0$ means that the trivial extension $B^c = B \ltimes I$
can be considered also as the tensor algebra of the $B$-$B$-bimodule $I$.

    \bigskip
Let us consider the matrix ring $B_2 = \bmatrix B & I \cr 0 & B\endbmatrix$.
The $B_2$-modules
can be written in the form 
$$
 (N_0,N_1;\gamma\:I\otimes_B N_1 \to N_0)
$$ 
where $N_0,N_1$ are
$B$-modules and $\gamma$ is a $B$-homomorphism; in terms
of matrices, we write $N$ in the form $\bmatrix N_0 \cr N_1\endbmatrix$ (and then
we can use matrix multiplication, taking into account the map $\gamma$).
   \bigskip
Consider the subring $B^2 = B\times B = \bmatrix B & 0 \cr 0 & B\endbmatrix$ of $B_2$.
We will identify $\mod B$ with the $B_2$-modules of the form $(N,0;0)$, thus with
those annihilated by the idempotent
$e_1 = \left[\smallmatrix 0 & 0 \cr 0 & 1\endsmallmatrix \right]$. 
We also consider a second embedding functor $\mod B \to \mod B_2$; it 
sends the $B$-module $N$ to 
$N[1] = (0,N;0)$; the $B_2$-modules of the form $N[1]$ are just the 
$B_2$-modules annihilated by the idempotent
$e_0 = \left[\smallmatrix 1 & 0 \cr 0 & 0\endsmallmatrix \right]$ (the reason for
writing $N[1]$ will become clear when we deal with $B_\infty$).
	\medskip
Let $\widetilde {\Cal N}$ denote the indecomposable $B_2$-modules 
which are not $B^2$-modules,
thus those indecomposable $B_2$-modules $N$ with $IN \neq 0$ (again, we take just
one module from each isomorphism class). 
    \medskip
{\bf Lemma 7.} {\it For any $B_2$-module $N = (N_0,N_1;\gamma)$ 
in $\widetilde {\Cal N}$, we have $N_0\in \add\Cal X$ and $N_1\in \add\Cal Y.$}
   \medskip
Note that this means: {\it There exists 
an exact sequence $0 \to X \to N \to Y[1] \to 0$ with $X\in \add\Cal X$ and
$Y \in \add\Cal Y,$} namely the sequence $(N_0,0;0) \to N \to (0;N_1;0).$
   \medskip
Proof: We have shown in Lemma 6 
that $I\otimes_B N_1$ belongs to $\add\Cal X$. 
Decompose $N_0 = X\oplus N'_0$ with $X$ in $\add\Cal X$
and $N'_0\in \add\Cal Y$. Since $\Hom(\Cal X,N'_0) = 0$, we 
see that we can split off $(N'_0,0;0)$ from $N$, thus $N'_0 = 0$ and therefore 
$N_0 = X$ belongs to $\add\Cal X$.

Instead of looking at the homomorphism
$\gamma\:I\otimes_B N_1 \to N_0$, we also may consider the adjoint map
$\gamma'\:N_1 \to \Hom_B(I,N_0).$ Since $\Hom_B(I,N_0)$ belongs to $\add\Cal Y$,
any direct summand of $N_1$ which belongs to $\add\Cal X$ has to lie in the kernel
of $\gamma'$, thus can be split off. This shows that $N_1$ belongs to $\add\Cal Y.$ 

   \bigskip
We may illustrate the structure of $\mod B_2$ in the following way:
$$
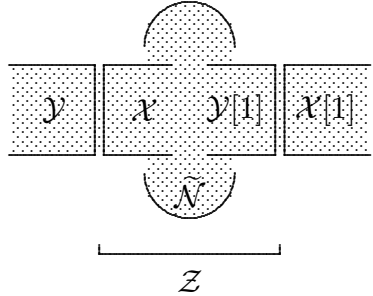
\hbox{\beginpicture
\setcoordinatesystem units <1cm,1cm>
\multiput{} at 0 -1  0 1 /
\put{$\Cal Y$} at  0 0
\put{$\Cal X$} at  1 0
\put{$\Cal Y[1]$} at 2 0
\put{$\Cal X[1]$} at 3 0
\put{$\widetilde {\Cal N}$} at 1.5 -.9
\plot -0.5 -0.5 0.45 -0.5 0.45 0.5  -0.5 0.5 /
\plot 1.3 -0.5 0.55 -0.5 0.55 0.5 1.3 0.5 /
\plot 3.5 -0.5 2.55 -0.5 2.55 0.5  3.5 0.5 /
\plot 1.7 -0.5 2.45 -0.5 2.45 0.5 1.7 0.5 /
\circulararc -180 degrees from 2 -.7 center at 1.5 -.7
\circulararc  180 degrees from 2  .7 center at 1.5 .7
\setshadegrid span <.5mm>
\vshade -.5 -.5 .5 <,z,,>
        1 -.5 .5  <z,z,,>
        1.05 -1 1    <z,z,,>
        1.2 -1.1 1.1    <z,z,,>
        1.3 -1.2 1.2    <z,z,,>
        1.5 -1.2 1.2    <z,z,,>
        1.7 -1.2 1.2    <z,z,,>
        1.8 -1.1 1.1    <z,z,,>
        1.95 -1 1    <z,z,,>
        2 -.5 .5    <z,,,>
        3.5 -.5 .5 /
\plot 0.5 -1.5  0.5 -1.6  2.5 -1.6  2.5 -1.5 /
\put{$\Cal Z$} at 1.5 -1.9 
\endpicture}
$$
Here, $\Cal Z$ denotes the indecomposable $B_2$-modules $N$
with an exact sequence 
$0 \to X \to N \to Y[1] \to 0$ with $X\in \add\Cal X$ and
$Y \in \add\Cal Y.$ Note that 
$\Cal Z$ consists of 
$\Cal X, \Cal Y[1]$ and $\widetilde{\Cal N}$. 
      \bigskip\bigskip
{\bf 6. The equivalence of $\Mat \Hom(\Cal G,\tau\Cal F)$ and $\Cal Z$.}
     \medskip
The second bimodule to be considered is $\Hom(\Cal G,\tau\Cal F).$ Here,
$\Cal A = \Cal G$, $\Cal B = \Cal F$ and $E(X,Y) = \Hom(X,\tau Y).$
      \medskip
{\bf Proposition 2.} {\it There is an equivalence of categories 
$$
 \eta\:\Mat \Hom(\Cal G,\tau\Cal F) \longrightarrow \Cal Z
$$
such that $\eta(M_1,M_2,\phi) = (G(M_1),F(M_2),\phi')$ for some $\phi'$.}
     \medskip
Proof: Let $(M_1,M_2,\phi)$ be given, with $M_1\in \Cal G,
M_2\in \Cal F, \phi\in \Hom(M_1,\tau M_2).$ Applying $G = \Hom_A(T,-)$ to $\phi$,
we obtain
$$
\align
 G(\phi)\:G(M_1) \to G(\tau M_2) &= \Hom_A(T,\tau M_2) \cr
 \cong &\Hom_A(\tau^{-1}T,M_2) \cr
 \cong &\Hom_B(F(\tau^{-1}T),F(M_2)) = \Hom_B(I,F(M_2)).
\endalign
$$
First, we have used again that $\tau^{-1}$ and $\tau$ are adjoint, and then the fact that
$F$ yields a bijection $\Hom_A(X,M_2) \to \Hom_B(F(X),F(M_2)),$ for any $A$-module $X$,
since $M_2$ belongs to $\Cal F$ (note that this bijection is $B$-linear, since 
$\tau^{-1}T$ is an $A$-$B^{\text{op}}$-bimodule).
The required map $\phi'$ is the adjoint of the map $G(M_1) \to \Hom_B(I,F(M_2)).$

For the converse, we only have to observe that $G$ yields a bijection
$\Hom_A(M_1,X) \to \Hom_B(G(M_1),G(X))$ for any $A$-module $X$, since $M_1$ belongs
to $\Cal G.$
   \bigskip\bigskip

{\bf 7. The algebras $B^c$ and $B_\infty$.}
     \medskip
We also consider the $(\Bbb Z\times \Bbb Z)$-matrix ring
$$
B_\infty = 
\bmatrix \ddots &  \ddots\cr
                & B & I \cr
                &   & B & I \cr
                &   &   & B      & \ddots \cr
                &   &   &  & \ddots
\endbmatrix
$$
of $(\Bbb Z\times \Bbb Z)$-matrices with only finitely many non-zero entries: 
on the main diagonal, there are copies of $B$, above the main diagonal, there
are copies of $I$ (since $I\otimes I = 0$, we do not have to worry about 
multiplying elements from different copies of $I$).
Note that here we deal 
with a ring without identity element, 
but at least it has sufficiently many idempotents. 
The $B_\infty$-modules are of the form $(N_i,\gamma_i)_i$,
indexed over $\Bbb Z$, with $B$-modules $N_i$ and $B$-linear maps 
$\gamma_i\:I\otimes_B N_i \to N_{i-1}$. Note that $B_\infty$ is locally
bounded (this means that any simple $B_\infty$-modules has a projective
cover and an injective envelope, both being of finite length).
      \medskip
We will consider $\mod B_2$ as the full subcategory of $\mod B_\infty$
with objects $(N_i,\gamma_i)$ where $N_i = 0$ for $i\notin\{0,1\}$.
In this way, we consider $\Cal Z$ as a fixed subcategory of $\mod B_\infty.$
Also, we define for any $t\in \Bbb Z$ a shift functor
$[t]\:\mod B_\infty \to \mod B_\infty$ by $N[t] = (N_{i-t},\gamma_{i-t})_i$,
for $N = (N_i,\gamma_i)_i$. 
    \medskip
{\bf Proposition 3.} {\it The indecomposable $B_\infty$-modules of finite
length are of the form
$N[t]$ with $N\in \Cal Z$ and $t\in \Bbb Z,$ and $N, t$ are uniquely determined.} 
       \medskip
Proof: Let $N = (N_i,\gamma_i)_i$ be a $B_\infty$-module of finite
length.  Decompose $N_i = X_i\oplus Y_i$ 
with $X_i\in \Cal X$ and $Y_i\in \Cal Y$, for all $i\in \Bbb Z$.
The discussion of the $B_2$-modules shows that $\gamma_i$ vanishes on $I\otimes X_i$
and maps into $Y_{i-1}$, thus we see that there is the following direct decomposition
of $B_\infty$-modules $M^{(t)} = (M^{(t)}_i,\gamma^{(t)}_i)_i$ where 
$$
  M^{(t)}_i = \left\{ \matrix \ X_{t-1} & & \quad\ \;i=t-1 \cr
                              Y_t &\text{for} & i = t \cr
                              0   &           & \qquad\text{otherwise}\endmatrix \right. 
$$
and where $\gamma^{(t)}_t$ is given by $\gamma_t.$ Obviously, $M^{(t)}$
belongs to $\Cal Z[t]$. 
As a consequence, if $N$ is an indecomposable $B_\infty$-module of finite length, 
then it belongs to $\Cal Z[t]$ for some $t\in \Bbb Z.$
     \bigskip
Recall that  a locally bounded ring $R$ is said to be {\it locally support-finite}
provided for any simple $R$-module $S$ there exists a finite set of simple
$R$-modules $\Cal S(S)$ with the following property: if $M$ is an
indecomposable $R$-module of finite length which has $S$ as a composition factor,
then any composition factor of $M$ belongs to $\Cal S(S)$.
     \medskip
{\bf Corollary.} {\it The algebra $B_\infty$ is locally support-finite.}
     \medskip
Proof. Let $S$ be a simple
$B_\infty$-module, say $S = (S_i,0)_i$ with $S_i = 0$ for $i\neq t$, and assume that
$N = (N_i,\gamma_i)_i$
is an indecomposable $B_\infty$-module of finite length with composition factor
$S$. Then $N$ belongs to $\Cal Z$ or to $\Cal Z[-1].$ This shows that there are
only finitely many simple $B^c$-modules which can occur as composition factors
of $N$. 

   \bigskip\bigskip
Our interest lies in the cluster-tilted algebra
$B^c = B \ltimes I$. Obviously, the algebra $B_\infty$ is a Galois covering
of $B^c$ with Galois group $\Bbb Z$ given by the shift functors $[t]$ with 
$t\in \Bbb Z$.
The covering functor $\pi\:\mod B_\infty \to \mod B^c$ sends
$(N_i,\gamma_i)_i$ to $(\bigoplus_i N_i,\gamma)$ with $\gamma$ being given by the
$\gamma_i$.
	\medskip
According to Dowbor-Lenzing-Skowronski [DLS], the Proposition 3 and its
Corollary have the following consequences:
	  \medskip
{\bf Corollary.} {\it The covering functor $\pi$ is dense and induces a
bijection between $\Cal Z$ and $\mod B^c$.}
	  \medskip
Of course, this bijection yields a bijection from $\widetilde{\Cal N}$ onto $\Cal N(B)$.
   \bigskip
Besides the covering functor $\pi\:\mod B_\infty \to \mod B^c$ itself, we also may
look at its restriction to $\mod B_2$.
Note that the subring of $B_2$ of all matrices of the form 
$\bmatrix b & x \cr 0 & b\endbmatrix$ with $b\in B$ and $x\in I$ is just 
$B^c$, and this inclusion gives rise to the restriction of the covering functor 
$$
 \mod B_2 \ \subset \  \mod B_\infty @>\pi>> \mod B^c.
$$

	\bigskip
The original definition of a cluster-tilted algebra $C$ as introduced by
Buan-Marsh-Reiten [BMR] implies that
the module category $\mod C$
is a factor category of a cluster category $\Cal C_A$. Namely, one
starts with the derived category $D^b(\mod A)$, say with shift functor $\Sigma$
and Auslander-Reiten translation $\tau$, and considers the orbit category
$\Cal C_A = D^b(\mod A)/\sigma$ with $\sigma = \Sigma\tau^{-1}.$ 
Then, one takes the factor category 
$\Cal C_A/\langle \tau T\rangle$ (here, $\langle \tau T\rangle$ is the ideal
of all maps which factor through $\add \tau T$).
It turns out that $T$, considered as an object of
$\Cal C_A/\langle \tau T\rangle$ is a progenerator, and its endomorphism ring
is $B^c$, thus one can identify
$$ 
 \mod B^c = \Cal C_A/\langle \tau T\rangle.
$$

We may change the procedure slightly: Starting with the derived category
$D^b(\mod A)$, we now first want to factor out the ideal 
$\langle \sigma^z(\tau T)\mid z\in \Bbb Z\rangle $ and only in the second step
form the orbit category with respect to the action of $\sigma$. We can 
make the identification
$$
 \mod B_\infty = D^b(\mod A)/\langle \sigma^z(\tau T)\mid z\in \Bbb Z\rangle
$$
so that the shift functor $M\mapsto M[1]$ 
on the left coincides with the operation of $\sigma$
on the right. In particular, the covering functor 
$$
 \pi\: \mod B_\infty \to \mod B_\infty/[1] = \mod B^c
$$ 
is nothing else than forming the orbit category 
$\left(D^b(\mod A)/\langle \sigma^z(\tau T)\mid z\in \Bbb Z\rangle\right)
 /\sigma.$
 \medskip
Let us remark that the importance of $B_2$ and $B_\infty$ for dealing with
a cluster-tilted algebra $B^c$ has been stressed already in [R].
  \bigskip\bigskip

{\bf 8. The category $\Mat \Hom(\Cal G,\tau\Cal F)$ for $T$ preprojective.}
     \medskip
In case $T$ is preprojective, we can improve the assertion of Lemma 7. 
   \medskip
{\bf Lemma 8.} {\it Assume that $T$ is preprojective.
For any $B_2$-module $N = (N_0,N_1,\gamma)$ in 
$\widetilde {\Cal N}$ we have $N_0\in \add \Cal X$ and 
$N_1 \in \add G(\Cal G\cap\Cal D)$.}
     \medskip
Proof: According to Lemma 7, we know that $N_0\in \add\Cal X$ and
$N_1\in \add \Cal Y$, thus $N_0 = F(M_0)$ for some $M_0\in \add \Cal F$
and $N_1 = G(M_1)$ for some $M_1\in \add\Cal G$. Instead of looking at
$\gamma\:I\otimes N_1 \to N_0$ we look again at the adjoint map
$$
 \gamma'\: N_1 = G(M_1) \longrightarrow \Hom_B(I,N_0) = \Hom_B(I,F(M_0)) 
$$
Using Lemma 5, we see that 
$$ 
 \Hom_B(I,F(M_0)) = \Hom_A(T,\tau M_0) = \Hom_A(T,t\tau M_0) = G(t\tau M_0).
$$ 
Since $G$
is an equivalence between $\Cal G$ and $\Cal Y$, there is a
homomorphism $f\:M_1 \to t\tau M_0$ such that $\gamma' = G(f).$ 

Now we use that $N$ belongs to $\widetilde{\Cal N}$.
If $M_1'$ is an indecomposable direct summand of $M_1$, then there must exist
an indecomposable direct summand $M_0'$ of $M_0$ such that 
$\Hom(M_1',t\tau M_0') \neq 0.$ Note that $M_0'$ belongs to $\Cal F$, thus 
to $\Cal D$, according to Lemma 2. But then also $M_1'$ belongs to $\Cal D$,
since by definition, $\Cal D$ is closed under predecessors. This shows that
$M_1'$ belongs to $\Cal G \cap \Cal D$, thus $N_1$ belongs to 
$\add G(\Cal G \cap \Cal D).$

      \bigskip
{\bf Corollary.} {\it Let $T$ be preprojective. Then 
$\Mat\Hom(\Cal G\cap \Cal D,\tau\Cal F)$
is directed.}
   \medskip
Proof: Under the equivalence $\eta$, the matrix category 
$\Mat\Hom(\Cal G\cap \Cal D,\tau\Cal F)$ is mapped to a subcategory 
$\Cal Z'$ of $\Cal Z$. We claim that $\Cal Z'$ is directed. This is clear in
case $A$ is representation finite, since in this case we deal with a subcategory
of a factor category of $D^b(\mod A)$, and $D^b(\mod A)$ is directed, thus also
$\Cal Z'$ is directed. Thus, we can assume that $A$ is representation infinite (and
connected). 
Let us recall the structure of the categories $\mod A$
and $D^b(\mod A)$. The category $\mod A$ decomposes into three parts:
the preprojectives $\Cal P$, the regular modules $\Cal R$ and the preinjectives
$\Cal Q$. Looking at $D^b(\mod A)$, the subcategories $\Sigma^z(\Cal Q)$
and $\Sigma^{z+1}(\Cal P)$ combine to form a transjective component, and any
finite subcategory of such a component is directed. But $\Cal Z'$ is a factor category
of a finite subcategory of the transjective component with the objects
$\Cal Q$ and $\Sigma(\Cal P)$, thus $\Cal Z'$ is directed.
      \bigskip\bigskip
{\bf 9. The bijection between $\Cal M(T)$ and $\Cal N(B)$.}
     \medskip

{\bf Proposition 4.} {\it Let $T$ be preprojective. 
There is a bijection $\iota\:\ind A \to \ind B^c,$
such that for $M\in \ind A$, the restriction of $\iota(M)$ to $B$ is 
$G(M)\oplus F(M).$}
	    \medskip
{\bf Remark.} Note that for any $A$-module $M$, we have
$$
 G(M) = G(tM) \quad\text{and}\quad F(M) = F(M/tM).
$$
Thus we could write $\iota(M) = G(tM)\oplus F(M/tM).$ This would stress that
we deal with the middle terms of the exact sequences
$$
\gather 
 0 \to tM \to M \to M/tM \to 0,\cr
 0 \to F(M/tM) \to \iota(M) \to G(tM) \to 0.
\endgather 
$$
	\medskip
Proof: For $M$ in $\Cal F$, let $\iota(M) = F(M);$
for $M$ in $\Cal G$, let $\iota(M) = G(M);$ thus, it remains to consider
$M$ in $\Cal M(T)$. 
    \medskip
We consider the categories $\Cal A = \Cal F$, $\Cal B = \Cal G\cap \Cal D$
and the bimodule $E = \Ext^1(\Cal F,\Cal G\cap \Cal D).$ 
The quadratic form $r_E$ is defined on 
$$
 K = K(\Cal F,\oplus) \oplus K(\Cal G\cap \Cal D,\oplus)
$$
We also may consider the bimodule $E' = \Hom(\Cal G\cap \Cal D,\tau \Cal F)$ with
quadratic form $r_{E'}$ on $K$. According to Auslander-Reiten, the bimodules
$E$ and $E'$ are dual to each other,  thus the quadratic forms $r_E$ and $r_{E'}$
coincide.
	
The indecomposable objects both in $\Cal M(T)$ and in $\Cal N(B)$
correspond bijectively to the positive non-simple roots of the quadratic 
form $r_E = r_{E'}$,
according to the theorem of de la Pe\~na and Simson quoted in section 3.
This completes the proof.
     \bigskip
     \bigskip
{\bf Remark.} As we see, a key ingredient of the proof is the 
duality of the bimodules $\Ext^1(\Cal F,\Cal G)$ and as $\Hom(\Cal G,\tau \Cal F)$,
which is one of the basic results of Auslander-Reiten theory, since we have to
deal with the matrices over these bimodules.  
(Note that we could write $\overline{\Hom}$ instead of $\Hom$,
since the only maps from injective modules to $\tau\Cal F$ are the zero maps.)
      \bigskip
{\bf Corollary.} {\it Let $T$ be preprojective. Let $Z \in \Cal N(B)$. Then
$\End(Z)= k$.}
     \medskip
Proof. We can write $Z = \pi(N)$ for some indecomposable $B_2$-module
$N = (N_1,N_2,\gamma)$. Now, $\End_B(N) = k$. Also, $\Hom(N_1,N_2) = 0$ 
according to the separation property. 
Thus $\End(Z) = k.$ 
     \medskip
Observe that this is the assertion (c) of Theorem 2.

	\bigskip\bigskip
{\bf 10. The quadratic form $q_B$.} 
     \medskip
Given a finite-dimensional algebra $R$ of finite global dimension, we denote
by $\langle -,-\rangle$ the bilinear form on $K_0(R)$ with
$$
 \langle \bdim M,\bdim M'\rangle = \sum_{t\ge 0} (-1)^t\dim_k\Ext_R^t(M,M')
$$
for all $R$-modules $M,M'$, and we write $q_B(x) = \langle x,x\rangle$ for $x\in
K_0(R)$; in this way, we obtain a quadratic form which is called the 
{\it Euler form.} 
     \medskip
Let us return to the hereditary algebra $A$ with tilting module $T$ and $B = \End(T)$.
Recall that we have denoted by $T_1,\dots,T_n$ indecomposable direct summands of $T$, 
one from each isomorphism class and   
$g\:K_0(A) \to K_0(B)$ was defined to be the linear bijection such that
$g(\bdim  T_i) = \bdim  G(T_i)$, for $1 \le i \le n$. 
	  \medskip
{\bf Addendum to Proposition 4.} {\it We have $\bdim \iota(M) = \abs g(\bdim M).$}
     \medskip
Proof. Since $tM$ belongs to $\Cal G$, we have $\bdim G(M) = \bdim G(tM) 
= g(\bdim tM)$. Since $M/tM$ belongs to $\Cal F$, we have $\bdim F(M) = \bdim F(M/tM) = 
- g(\bdim M/tM)$. Altogether, it follows from $\bdim M = \bdim tM +\bdim M/tM$ 
that
$$ 
 g(\bdim M) = g(\bdim tM) +g(\bdim M/tM) = g(\bdim G(M)) - g(\bdim F(M)).
$$ 
The separation property now implies that 
$$ 
 \abs g(\bdim M) =  g(\bdim G(M)) + g(\bdim F(M)).
$$ 
   \medskip

{\bf Lemma 9.} Let $M \in \Cal M(T)$. Then 
$$
 q_A(\bdim tM) = \dim \End (tM),\quad q_A(\bdim M/tM) = \dim \End(M/tM).
$$
	\medskip
Proof: Since $A$ is hereditary, $q_A(\bdim X) = \dim \End(X) - \dim\Ext^1(X,X)$
for any $A$-module $X$. According to Lemma 4, both $tM$ and $M/tM$ are
modules without self-extensions. 
	\bigskip
We have denoted by 
$T_1,\dots,T_n$ indecomposable direct summands of $T$, one from each
isomorphism class. If we define 
$$
 g\:K_0(A) \to K_0(B) \quad\text{by}\quad  g(x) = (\langle t_i,x\rangle)_i,
$$
then $g(\bdim  T_i) = \bdim  G(T_i)$, for $1 \le i \le n$, thus $g$ is the linear
bijection between $K_0(A)$ and $K_0(B)$ mentioned in the introduction.
	  \medskip

{\bf Proposition 5.} {\it Let $M \in \Cal M(T)$.}	
$$
 q_B(\abs g(\bdim M)) = 2(\dim\End(tM)+\dim\End(M/tM))-1.
$$
	\medskip
Proof: Let  $x = \dim M.$ Since $M$ is preprojective, 
$q_A(x) = 1.$ Write $y = \bdim tM,$ and $z = \bdim M/tM$.
Since $x = y+z$ and $g$ is linear, $g(x) = g(y)+g(z)$.
Now $g(y)$ is positive and $g(z)$ is negative. Since the support of $F(M)$ and
$G(M)$ are disjoint (the separation property), we see that
$\abs g(x) = g(y) - g(z).$ Thus
$$
\align
 q_B(\abs g(x)) &= q_B(g(y)-g(z)) \cr
    &= q_B(g(y)) + q_B(g(z)) - 2(g(y),g(z))_B \cr
    &= q_A(y) + q_A(z) - 2(y,z)_A.
\endalign
$$
On the other hand,
$$
 1 = q_A(x) = q_A(y+z) =  q_A(y) + q_A(z) + 2(y,z)_A.
$$
If we add the two equalities, we obtain
$$
 q_B(\abs g(x)) + 1 = 2( q_A(y) + q_A(z)).
$$
Now, we apply Lemma 9 in order to see that
$$
  q_B(\abs g(x)) + 1 = 2(\dim\End(tM)+\dim\End(M/tM)).
$$
This completes the proof.

     \medskip
{\bf Corollary.} {\it If $M\in \Cal M(T)$, then $q_B(\abs g(\bdim M))$ is an
odd integer greater than $2$.}
    \medskip
This shows the assertion (b) of Theorem 2.
     \medskip 
Proof: For $M\in \Cal M(T)$, both modules $tM$ and $M/tM$ are non-zero, and therefore
$\dim\End(tM) \ge 1$ and $\dim\End(M/tM) \ge 1.$ 
	      \medskip
{\bf Corollary.} {\it If $M\in \Cal M(T)$, then $q_B(\abs g(\bdim M)) = 3$
if and only if $tM$ and $M/tM$ are both indecomposable.}
   \medskip
We will discuss this situation in the next section.
   \bigskip\bigskip
{\bf 11. The mixed pairs.}
     \medskip
In order to determine the quadratic form $r_E$, one needs to know the pairs
$(X,Y)$ of indecomposable $A$-modules with $X\in \Cal F$ and $Y\in \Cal G$
such that $\Ext^1(X,Y) \neq 0.$ 
     \medskip
We call $(X,Y)$ a {\it mixed pair} provided $X$ is an indecomposable 
module in $\Cal F(T)$, whereas
$Y$ is an indecomposable module in $\Cal G(T)$, and finally $\Ext^1(X,Y) \neq 0.$ 
    \medskip

{\bf Proposition 6.} {\it Let $T$ be preprojective. Then any 
mixed pair $(X,Y)$ is an orthogonal exceptional pair consisting of
preprojective $A$-modules such that  $\Ext^1(X,Y)$ is
one-dimensional. Let $M$ be the middle term of a non-split exact 
sequence $0 \to Y \to M \to X \to 0$, then $M$
is indecomposable and preprojective, and 
$\Ext^1(X,M) = 0$ and  $\Ext^1(M,Y) = 0.$}
	     \medskip
Proof: Let $0 \to Y \to M \to X \to 0$ be a non-split exact sequence. Since
$\Hom(Y,X) = 0$, it follows that $M$ is indecomposable. Since $M$ is
a predecessor of $X$, we see that $M$ is preprojective. Since $Y$ is a proper 
predecessor of $X$, it follows that $\Hom(X,Y) = 0$ and $\Ext^1(Y,X) = 0$, thus
$(X,Y)$ is an orthogonal and exceptional pair. The 
full subcategory $\Cal C$ of modules with
a filtration with factors $X$ and $Y$ consists of $X$, $Y$ and some modules in
$\Cal M(T)$,
thus it is of finite type, therefore $\Ext^1(X,Y)$ is at most one-dimensional,
thus one-dimensional.

It is well-known (and easy to see) that $\Cal C$ is equivalent to the category 
of representations of the quiver of Dynkin type $\Bbb A_2$, and
$M$, considered as an object of $\Cal C$ is both projective and injective.
In particular, we have 
$\Ext^1(X,M) = 0$ and  $\Ext^1(M,Y) = 0.$
	     \bigskip
{\bf Remark.} There is the following converse: {\it If $(X,Y)$ is an orthogonal
exceptional pair consisting of preprojective $A$-modules, then the dimension
of $\Ext^1(X,Y)$ is at most $1$. If $\dim_k\Ext^1(X,Y) = 1$, then there is a
preprojective tilting module $T$ with $X\in \Cal F(T)$ and $Y\in \Cal G(T)$} (and therefore
$(X,Y)$ is a mixed pair for $T$).
	\medskip
Proof. First, if $\dim_k\Ext^1(X,Y)\ge 2,$ then there are infinitely many indecomposable
$A$-modules $M$ which have a submodule $M'$ which is a direct sum of copies of $Y$
such that $M/M'$ is a direct sum of copies of $X$, and all these modules are predecessors
of $X$, this is impossible. This shows that  $\dim_k\Ext^1(X,Y)\le 1.$ 

Second, assume that $\dim_k\Ext^1(X,Y) = 1$.   
We claim that $\tau^{-1}X\oplus Y$ is a partial tilting module. Since $Y$ is a 
predecessor of $X$, it is also a predecessor of $\tau^{-1}X$, therefore
$\Ext^1(Y,\tau^{-1}X) = 0$. On the other hand, also
$\Ext^1(\tau^{-1}X,Y) = D\Hom(Y,X) = 0$. The Bongartz completion $T$ of this partial
tilting module is a preprojective tilting module and 
$X\in \Cal F(T)$ and $Y\in \Cal G(T).$

      \bigskip\bigskip
{\bf 12. Proof of Theorem 1.}
     \medskip
A quadratic form $q$ defined on $\Bbb Z^n$ with values in $\Bbb Z$
is said to be an {\it integral} form.
   \medskip
{\bf Proposition 7.}  {\it Let $q$ be an integral quadratic form on $\Bbb Z^n$
which is positive definite. If $x, x'\in \Bbb Z^n$ satisfy $q(x) = 1 = q(x')$ and
$\abs x = \abs x'$, then $x = \pm x'.$}
      \medskip
Proof. Let $y\in \Bbb Z^n$ be defined by $y_i = x_i$ provided $x_i = x'_i$ and
$y_i = 0$ otherwise. Let $z = x - y,$ thus $x = y+z$ and $x' = y - z$.
Let $(-,-)$ be the  bilinear form (with values in $\frac12 \Bbb Z$) corresponding to $q$.
Then
$$
\align
 1 &= q(x) = q(y+z) = q(y)+q(z) + 2(y,z),\cr
 1 &= q(x') = q(y-z) = q(y)+q(z) - 2(y,z).
\endalign
$$
shows that $(y,z) = 0,$ thus 
$$
 1 = q(x) = q(y)+q(z).
$$
Since we assume that $q$ takes values in $\Bbb Z$ and since 
$q$ is positive definite, it follows that $y = 0$ or $z = 0.$
If $z = 0$, then $x = x'$. If $y = 0,$ then $x = -x'$
   \medskip
Proof of Theorem 1. Let $C$ be representation-finite and cluster-tilted, say
$C = B^c$ with $B = \End_A(T)$, where $T$ is a tilting $A$-module and
$A$ is hereditary and representation-finite. Let $N,N'$ be indecomposable $C$-modules
with $\bdim N = \bdim N'.$ Proposition 4 provides indecomposable $A$-modules
$M, M'$ such that the restriction of $\iota(M)$ to $B$ is $N,$ and 
the restriction of $\iota(M')$ to $B$ is $N'.$ Let $x = \bdim M, x' = \bdim M'$.
Then $\bdim N = \abs g(x),$ and $\bdim N' = \abs g(x')$ according to the addendum
to Proposition 4. Now, $q_B(g(x)) = q_A(x) = 1,$
and also $q_B(g(x')) = q_A(x') = 1.$ With $q_A$ also $q_B$ is positive definite.
Thus we can apply Proposition 7 in order to see that $g(x) = \pm g(x')$ and
therefore $x = \pm x'.$ However, both $x, x'$ are non-negative vectors, thus $x = x'$
and therefore $M, M'$ are isomorphic (since any real root module is determined by 
its dimension vector).
Since $\iota$ is a bijection of isomorphism
classes, it follows that $N, N'$ are isomorphic. 
	 \bigskip\bigskip

{\bf 13. Examples.}
     \medskip
{\bf 13.1.}
Let us exhibit one example in detail. In particular, we will
see that the categories $\Cal M(T)$ and $\Cal N(B)$ can be quite different!

Consider the algebra $A = \Bbb T_{33}$; this is the
path algebra of the quiver $Q$ of type $\Bbb A_5$, with a unique sink and indecomposable
projective modules of length at most $3$. We label the vertices as exhibited on the left.
To the right, we present the Auslander-Reiten 
quiver and mark a tilting module using $*$: it consists of the indecomposable projective modules
of length $1$ and $3$ and the simple injective modules:
$$
 T(a) = P(1),\quad T(b) = P(3), \quad T(b') = P(3'),\quad T(c) = I(3), \quad T(c') = I(3').
$$ 
The class $\Cal G$ of indecomposable torsion modules consists of the modules $T(a), T(b), T(b')$ 
and the five indecomposable injective modules, the class $\Cal F$ of indecomposable 
torsionfree modules consists of the
two modules $\tau I(3)$ and $\tau I(3').$
$$
\hbox{\beginpicture
\setcoordinatesystem units <.7cm,.7cm>
\put{\beginpicture
\multiput{$\circ$} at 0 2  1 1  2 0  1 3  2 4 /
\arr{1.8 0.2}{1.2 0.8} 
\arr{1.8 3.8}{1.2 3.2} 
\arr{0.8 1.2}{0.2 1.8} 
\arr{0.8 2.8}{0.2 2.2} 
\put{$1$} at -0.5 2 
\put{$2$} at .5 3.2
\put{$2'$} at .5 0.8
\put{$3$} at 1.5 4.2 
\put{$3'$} at 1.5 -0.2
\put{$A$} at -0.3 4.5

\endpicture} at 0 0
\put{\beginpicture
\multiput{$\ast$} at 0 2  2 0  2 4  6 0  6 4 /
\multiput{$\bigcirc$} at 0 2  2 0  2 4  6 0  6 4 /
\plot 2 0  0 2  2 4 /
\plot 4 0  2 2  4 4 /
\plot 6 0  4 2  6 4 /
\plot 2 0  6 4 /
\plot 2 4  6 0 /
\plot 1 1  4 4  5 3 /
\plot 1 3  4 0  5 1 /
\multiput{$\Cal M$} at 1 1  1 3  2 2  3 1  3 3 /
\put{$T(a)$} at -.8 2
\put{$T(b)$} at 1 4
\put{$T(b')$} at 1 0
\put{$T(c)$} at 7 0
\put{$T(c')$} at 7 4
\setshadegrid span <.5mm>
\hshade 0 2 6 <,,,z> 2 0 4 <,,z,> 4 2 6 /

\endpicture} at 8 0
\put{} at -3 0
\endpicture}
$$
The positions of the five mixed modules are denoted by $\Cal M$, in a later table we
will provide more information on these modules.

Here is the quiver with relations for the
algebra $B = \End({}_AT)$, as well as
the Auslander-Reiten quiver
of $\mod B$. 
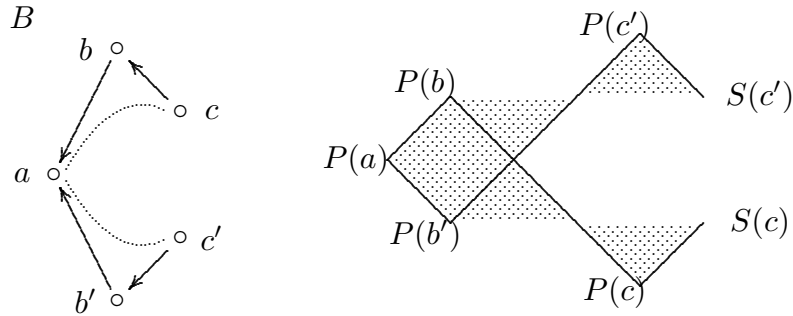
$$
\hbox{\beginpicture
\setcoordinatesystem units <.7cm,.7cm>
\put{\beginpicture
\multiput{$\circ$} at 0 2  1 0  1 4  2 1  2 3  /
\arr{0.9 0.2}{0.1 1.8} 
\arr{0.9 3.8}{0.1 2.2} 
\arr{1.8 0.8}{1.2 0.2} 
\arr{1.8 3.2}{1.2 3.8} 
\setdots <.5mm>
\put{$a$} at -.5 2
\put{$b$} at .5 4
\put{$b'$} at .5 0
\put{$c$} at 2.5 3
\put{$c'$} at 2.5 1
\put{$B$} at -0.5 4.5
\setquadratic
\setdots <.5mm>
\plot 0.2 2.1  1 3  1.8 3 /
\plot 0.2 1.9  1 1  1.8 1 /

\endpicture} at 0 0

\put{\beginpicture
\multiput{} at 0 0  6 4 /
\plot 5 1  4 0  1 3  0 2  1 1  4 4  5 3 /

\setshadegrid span <.5mm>
\hshade 0 4 4  1 3 5 /
\hshade 1 1 3 <,,,z> 2 0 2 <,,z,> 3 1 3 /
\hshade 3 3 5  4 4 4 /
\put{$P(a)$} at -0.5 2
\put{$P(b')$} at 0.6 0.8
\put{$P(b)$} at 0.6 3.2
\put{$P(c)$} at 3.6 -0.1
\put{$P(c')$} at 3.6 4.1

\put{$S(c')$} at 5.9 3
\put{$S(c)$} at 5.9 1

\endpicture} at 7 0
\endpicture}
$$
The $A$-modules in $\Cal G$ and $\Cal F$, and the corresponding $B$-modules 
obtained by applying the functors $G$ and $F$, respectively, are as follows:
$$
\matrix
P(1) & P(3) & P(3') & I(1) & I(2) & I(2') & I(3) & I(3') & & \tau I(3) & \tau I(3') \cr
P(a) & P(b) & P(b') & I(a) & S(b) & S(b')  & P(c) & P(c') & & S(c) & S(c')
\endmatrix
$$

Next, we present  the shape of the cluster category with circles showing the
direct summands of $T$, or better, just their labels 
(these modules are now considered as objects of $\Cal C_A$), 
always, the dashed lines have to be identified in order to form 
Moebius strips:
$$
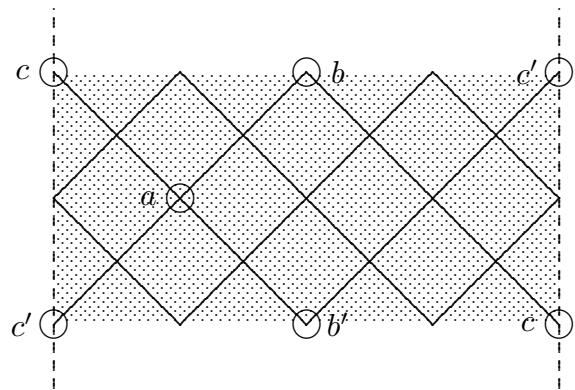
\hbox{\beginpicture
\setcoordinatesystem units <.7cm,.7cm>
\plot 0 2  2 0  6 4  8 2  6 0  2 4  0 2 /
\multiput{$\bigcirc$} at 0 0  4 0  8 0  0 4  4 4  8 4  2 2 /
\plot 0 0  4 4  8 0 /
\plot 0 4  4 0  8 4 /
\setdashes <1mm>
\plot 0 -1  0 5 /
\plot 8 -1  8 5 /
\setshadegrid span <.5mm>
\hshade 0 0 8  4 0 8 /
\put{$a$} at 1.5 2
\put{$b$} at 4.5 4
\put{$b'$} at 4.5 0
\put{$c$} at 7.5 0
\put{$c'$} at 7.5 4
\put{$c'$} at -.5 0
\put{$c$} at -.5 4
\multiput{} at 0 -0.5  0 4.5 /
\endpicture}
$$

Now, we present first the quiver with relations for the
algebra  $B^c$, as well as
the Auslander-Reiten quiver
of $\mod B^c$. 
The positions of the five mixed modules are denoted by $\Cal N$.
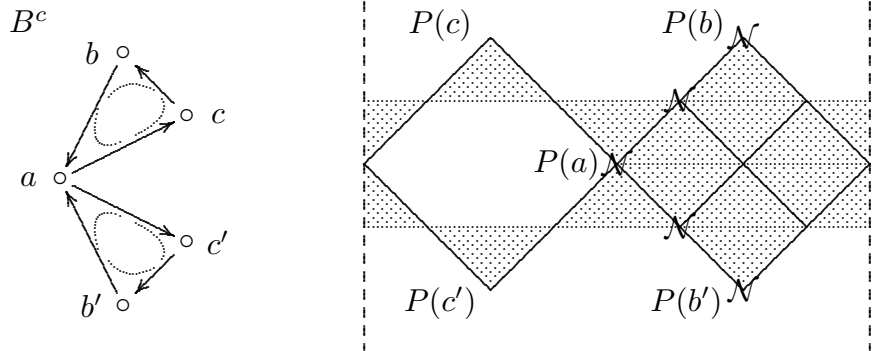
$$
\hbox{\beginpicture
\setcoordinatesystem units <.7cm,.7cm>

\put{\beginpicture
\multiput{$\circ$} at 0 2  1 0  1 4  2 1  2 3  /
\arr{0.9 0.2}{0.1 1.8} 
\arr{0.9 3.8}{0.1 2.2} 
\arr{1.8 0.8}{1.2 0.2} 
\arr{1.8 3.2}{1.2 3.8} 

\arr{0.2 2.1}{1.8 2.9}
\arr{0.2 1.9}{1.8 1.1}

\put{$a$} at -.5 2
\put{$b$} at .5 4
\put{$b'$} at .5 0
\put{$c$} at 2.5 3
\put{$c'$} at 2.5 1
\put{$B^c$} at -.5 4.5

\setquadratic
\setdots <.4mm>
\plot 0.77 3.3  1.1  3.5  1.5 3.4 /
\plot 0.6 3.1  0.6 2.6  1 2.6 /
\plot 1.2 2.7  1.6 3  1.6 3.3 /

\plot 0.77 0.7  1.1  0.5  1.5 0.6 /
\plot 0.6 0.9  0.6 1.4  1 1.4 /
\plot 1.2 1.3  1.6 1  1.6 0.7 /

\endpicture} at 5 0

\put{\beginpicture
\plot 0 2  2 0  6 4  8 2  6 0  2 4  0 2 /
\plot 5 1  7 3 /
\plot 5 3  7 1 /
\multiput{$\Cal N$} at 4 2  5 1  5 3  6 0  6 4 /
\setdashes <1mm>
\plot 0 -1  0 5 /
\plot 8 -1  8 5 /
\setdots <1mm>
\setdots <.5mm>

\plot 0 1  8 1 /
\plot 0 3  8 3 /
\plot 4 2  8 2 /
\setshadegrid span <.5mm>
\hshade  1 0 1  2 0 0  3 0 1 /
\hshade 0 2 2  1 1 3 /
\hshade 3 1 3  4 2 2 /
\hshade 1 3 5  2 4 4  3 3 5 /
\hshade 0 6 6  1 5 7  2 4 8  3 5 7  4 6 6 /
\hshade 1 7 8  2 8 8  3 7 8 /
\multiput{} at 0 -.5  0 4.5 /

\put{$P(a)$} at 3.2 2
\put{$P(b)$} at 5.2 4.2
\put{$P(b')$} at 5.1 -.2
\put{$P(c')$} at 1.2 -.2
\put{$P(c)$} at 1.2 4.2

\endpicture} at 12 0
\endpicture}
$$

	\bigskip	

The following table shows the bijection $\iota$ between
the modules $M$ in $\Cal M(T)$ and the modules $N = \iota(M)$ in 
$\Cal N(B)$. Below any $M$ we outline its torsion part $tM$ and
its torsionfree part $M/tM$ by writing
 $\dfrac {M/tM}{tM}$; similarly, under $N$ we show 
$\dfrac {N/tN}{tN}$. In the lowest row, one finds the values
$q_B(\bdim N).$ 
$$
\hbox{\beginpicture
\setcoordinatesystem units <1.2cm,2cm>

\put{$M$} at -1.5 7

\put{\beginpicture
\setcoordinatesystem units <.4cm,.4cm>
\multiput{$\circ$} at   2 0   2 4 /

\arr{0.7 1.3}{0.3 1.7} 
\arr{0.7 2.7}{0.3 2.3} 
\setdots <1mm>
\arr{1.7 0.3}{1.3 0.7} 
\arr{1.7 3.7}{1.3 3.3} 

\put{$1$} at 0 2
\put{$2$} at 1 3
\put{$2'$} at 1 1
\endpicture} at 0 7
\put{\beginpicture
\setcoordinatesystem units <.4cm,.4cm>
\multiput{$\circ$} at  1 1  2 0   2 4 /

\arr{0.7 2.7}{0.3 2.3} 
\setdots <1mm>
\arr{0.7 1.3}{0.3 1.7} 
\arr{1.7 0.3}{1.3 0.7} 
\arr{1.7 3.7}{1.3 3.3} 

\put{$1$} at 0 2
\put{$2$} at 1 3
\endpicture} at 2 7
\put{\beginpicture
\setcoordinatesystem units <.4cm,.4cm>
\multiput{$\circ$} at   2 0  1 3  2 4 /
\arr{0.7 1.3}{0.3 1.7} 

\setdots <1mm>
\arr{1.7 0.3}{1.3 0.7} 
\arr{1.7 3.7}{1.3 3.3} 
\arr{0.7 2.7}{0.3 2.3} 

\put{$1$} at 0 2
\put{$2'$} at 1 1
\endpicture} at 4 7
\put{\beginpicture
\setcoordinatesystem units <.4cm,.4cm>
\multiput{$\circ$} at   2 4 /

\arr{1.7 0.3}{1.3 0.7} 
\arr{0.7 1.3}{0.3 1.7} 
\arr{0.7 2.7}{0.3 2.3} 
\setdots <1mm>
\arr{1.7 3.7}{1.3 3.3} 

\put{$1$} at 0 2
\put{$2$} at 1 3
\put{$2'$} at 1 1
\put{$3'$} at 2 0
\endpicture} at 6 7
\put{\beginpicture
\setcoordinatesystem units <.4cm,.4cm>
\multiput{$\circ$} at   2 0   /

\arr{1.7 3.7}{1.3 3.3} 
\arr{0.7 1.3}{0.3 1.7} 
\arr{0.7 2.7}{0.3 2.3} 
\setdots <1mm>
\arr{1.7 0.3}{1.3 0.7} 

\put{$1$} at 0 2
\put{$2$} at 1 3
\put{$2'$} at 1 1
\put{$3$} at 2 4
\endpicture} at 8 7

\put{\beginpicture
\setcoordinatesystem units <.4cm,.4cm>
\put{$\dfrac {\tau I(3)\oplus\tau I(3')}{P(1)}$} at 0 0
\endpicture} at 0 6

\put{\beginpicture
\setcoordinatesystem units <.4cm,.4cm>
\put{$\dfrac {\tau I(3)}{P(1)}$} at 0 0
\endpicture} at 2 6

\put{\beginpicture
\setcoordinatesystem units <.4cm,.4cm>
\put{$\dfrac {\tau I(3')}{P(1)}$} at 0 0
\endpicture} at 4 6

\put{\beginpicture
\setcoordinatesystem units <.4cm,.4cm>
\put{$\dfrac {\tau I(3)}{P(3')}$} at 0 0
\endpicture} at 6 6

\put{\beginpicture
\setcoordinatesystem units <.4cm,.4cm>
\put{$\dfrac {\tau I(3')}{P(3)}$} at 0 0
\endpicture} at 8 6

\put{$N$} at -1.5 5

\put{\beginpicture
\setcoordinatesystem units <.4cm,.4cm>
\arr{0.4 2.2}{1.6 2.8}
\arr{0.4 1.8}{1.6 1.2}
\setdots <1mm>
\arr{0.8 0.4}{0.2 1.6} 
\arr{0.8 3.6}{0.2 2.4} 
\arr{1.7 0.7}{1.3 0.3} 
\arr{1.7 3.3}{1.3 3.7} 
\put{$a$} at 0 2
\put{$\circ$} at 1 4
\put{$\circ$} at 1 0
\put{$c$} at 2 3
\put{$c'$} at 2 1

\endpicture} at 0 5

\put{\beginpicture
\setcoordinatesystem units <.4cm,.4cm>
\multiput{$\circ$} at 1 0  1 4  2 1    /
\arr{0.4 2.2}{1.6 2.8}
\setdots <1mm>
\arr{0.8 3.6}{0.2 2.4} 
\arr{0.8 0.4}{0.2 1.6} 
\arr{1.7 0.7}{1.3 0.3} 
\arr{1.7 3.3}{1.3 3.7} 
\arr{0.4 1.8}{1.6 1.2}
\put{$a$} at 0 2
\put{$c$} at 2 3
\endpicture} at 2 5

\put{\beginpicture
\setcoordinatesystem units <.4cm,.4cm>
\arr{0.4 1.8}{1.6 1.2}
\setdots <1mm>
\arr{0.4 2.2}{1.6 2.8}
\arr{0.8 0.4}{0.2 1.6} 
\arr{0.8 3.6}{0.2 2.4} 
\arr{1.7 0.7}{1.3 0.3} 
\arr{1.7 3.3}{1.3 3.7} 
\put{$a$} at 0 2
\put{$\circ$} at 1 4
\put{$\circ$} at 1 0
\put{$\circ$} at 2 3
\put{$c'$} at 2 1
\endpicture} at 4 5

\put{\beginpicture
\setcoordinatesystem units <.4cm,.4cm>

\multiput{$\circ$} at   1 4  2 1   /
\arr{0.4 2.2}{1.6 2.8}
\arr{0.8 0.4}{0.2 1.6} 
\setdots <1mm>
\arr{0.8 3.6}{0.2 2.4} 
\arr{1.7 0.7}{1.3 0.3} 
\arr{1.7 3.3}{1.3 3.7} 
\arr{0.4 1.8}{1.6 1.2}
\put{$a$} at 0 2
\put{$b'$} at 1 0
\put{$c$} at 2 3
\endpicture} at 6 5

\put{\beginpicture
\setcoordinatesystem units <.4cm,.4cm>
\multiput{$\circ$} at   1 0    2 3  /
\arr{0.8 3.6}{0.2 2.4} 
\arr{0.4 1.8}{1.6 1.2}
\setdots <1mm>
\arr{0.8 0.4}{0.2 1.6} 
\arr{1.7 0.7}{1.3 0.3} 
\arr{1.7 3.3}{1.3 3.7} 
\arr{0.4 2.2}{1.6 2.8}
\put{$a$} at 0 2
\put{$b$} at 1 4
\put{$c'$} at 2 1
\endpicture} at 8 5

\put{\beginpicture
\setcoordinatesystem units <.4cm,.4cm>
\put{$\dfrac {GP(1)}{F\tau I(3)\oplus F\tau I(3')}$} at 0 0
\endpicture} at 0 4

\put{\beginpicture
\setcoordinatesystem units <.4cm,.4cm>
\put{$\dfrac {GP(1)}{F\tau I(3)}$} at 0 0
\endpicture} at 2 4

\put{\beginpicture
\setcoordinatesystem units <.4cm,.4cm>
\put{$\dfrac {GP(1)}{F\tau I(3')}$} at 0 0
\endpicture} at 4 4

\put{\beginpicture
\setcoordinatesystem units <.4cm,.4cm>
\put{$\dfrac {GP(3')}{F\tau I(3)}$} at 0 0
\endpicture} at 6 4

\put{\beginpicture
\setcoordinatesystem units <.4cm,.4cm>
\put{$\dfrac {GP(3)}{F\tau I(3')}$} at 0 0
\endpicture} at 8 4

\put{$q_B(\bdim N)$} at -1.5 3.2
\put{\beginpicture
\setcoordinatesystem units <.4cm,.4cm>
\put{$5$} at 0 0
\endpicture} at 0 3.2

\put{\beginpicture
\setcoordinatesystem units <.4cm,.4cm>
\put{$3$} at 0 0
\endpicture} at 2 3.2

\put{\beginpicture
\setcoordinatesystem units <.4cm,.4cm>
\put{$3$} at 0 0
\endpicture} at 4 3.2

\put{\beginpicture
\setcoordinatesystem units <.4cm,.4cm>
\put{$3$} at 0 0
\endpicture} at 6 3.2

\put{\beginpicture
\setcoordinatesystem units <.4cm,.4cm>
\put{$3$} at 0 0
\endpicture} at 8 3.2

\endpicture}
$$

	\bigskip

Finally, we note that the quadratic form $r_E$ has the following graph:
$$
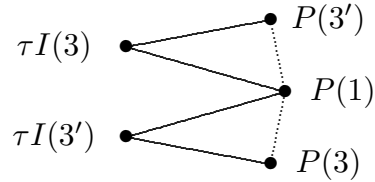
\hbox{\beginpicture
\setcoordinatesystem units <.8cm,1cm>
\multiput{$\bullet$} at 0 0  0 1  2 -0.3  2.2 0.5  2 1.3 /
\plot 0 0  2 -0.3 /
\plot 0 0  2.2  0.5 /
\plot 0 1  2.2  0.5 /
\plot 0 1  2  1.3 /
\setdots <.5mm>
\plot 2 -0.3 2.2 0.5 /
\plot 2  1.3 2.2 0.5 /
\put{$\tau I(3)$} at -1 1
\put{$\tau I(3')$} at -1 0
\put{$P(3')$} at 2.8 1.3
\put{$P(1)$} at 3 .5
\put{$P(3)$} at 2.8 -0.3
\endpicture}
$$
As usual, we have deleted the isolated vertices (here, a vertex 
is said to be {\it isolated} provided it is not the endpoint of any solid
edge).
	\bigskip
	\bigskip
{\bf 13.2.} Next, let us present two non-isomorphic tilted algebras $B, B'$
such that the cluster-tilted algebras $B'$ and $(B')^c$ are isomorphic and
representation-finite.
$$
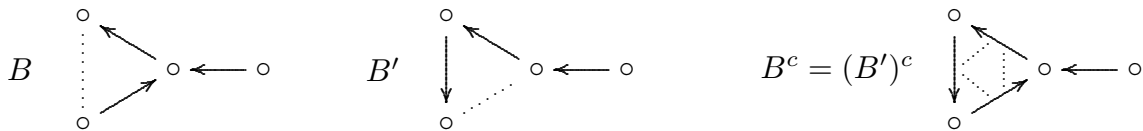
\hbox{\beginpicture
\setcoordinatesystem units <1cm,.6cm>
\put{\beginpicture
\put{$B$} at -.7 0
\multiput{$\circ$} at 0 1  0 -1  1 0  2 0 /
\arr{0.8 0.2}{0.2 0.8}
\arr{0.2 -0.8}{0.8 -0.2}
\arr{1.8 0}{1.2 0}
\setdots <1mm>
\plot 0 -0.7  0 0.7 /
\endpicture} at 0 0
\put{\beginpicture
\put{$B'$} at -.7 0
\multiput{$\circ$} at 0 1  0 -1  1 0  2 0 /
\arr{0.8 0.2}{0.2 0.8}
\arr{0 0.6}{0 -0.6} 
\arr{1.8 0}{1.2 0}
\setdots <1mm>
\plot 0.2 -0.8  0.8 -0.2 /
\endpicture} at 4 0
\put{\beginpicture
\put{$B^c = (B')^c$} at -1.3 0
\multiput{$\circ$} at 0 1  0 -1  1 0  2 0 /
\arr{0.8 0.2}{0.2 0.8}
\arr{0.2 -0.8}{0.8 -0.2}
\arr{0 0.6}{0 -0.6} 
\arr{1.8 0}{1.2 0}
\setdots <.7mm>
\plot 0.55 -0.3  0.55 0.3 /
\plot 0.1 0.1  0.4 0.5 /
\plot 0.1 -.1  0.4 -.5 /
\endpicture} at 9 0
\endpicture}
$$
There are 10 isomorphism classes of indecomposable $B^c$-modules $N$.
The following table lists the values of $q_B(\bdim N)$ and
$q_{B'}(\bdim N)$ for these modules.
$$
\hbox{\beginpicture
\setcoordinatesystem units <1cm,.7cm>
\put{\beginpicture
\setcoordinatesystem units <.15cm,.1cm>
\put{$\bdim N$} at 0 0
\endpicture} at -1.5 2
\put{\beginpicture
\setcoordinatesystem units <.15cm,.1cm>
\put{$\sssize 1$} at 0 1
\put{$\sssize 0$} at 0 -1
\put{$\sssize 0$} at 1 0
\put{$\sssize 0$} at 2 0
\endpicture} at 0 2
\put{\beginpicture
\setcoordinatesystem units <.15cm,.1cm>
\put{$\sssize 0$} at 0 1
\put{$\sssize 1$} at 0 -1
\put{$\sssize 0$} at 1 0
\put{$\sssize 0$} at 2 0
\endpicture} at 1 2
\put{\beginpicture
\setcoordinatesystem units <.15cm,.1cm>
\put{$\sssize 0$} at 0 1
\put{$\sssize 0$} at 0 -1
\put{$\sssize 1$} at 1 0
\put{$\sssize 0$} at 2 0
\endpicture} at 2 2
\put{\beginpicture
\setcoordinatesystem units <.15cm,.1cm>
\put{$\sssize 0$} at 0 1
\put{$\sssize 0$} at 0 -1
\put{$\sssize 0$} at 1 0
\put{$\sssize 1$} at 2 0
\endpicture} at 3 2
\put{\beginpicture
\setcoordinatesystem units <.15cm,.1cm>
\put{$\sssize 1$} at 0 1
\put{$\sssize 1$} at 0 -1
\put{$\sssize 0$} at 1 0
\put{$\sssize 0$} at 2 0
\endpicture} at 4 2
\put{\beginpicture
\setcoordinatesystem units <.15cm,.1cm>
\put{$\sssize 1$} at 0 1
\put{$\sssize 0$} at 0 -1
\put{$\sssize 1$} at 1 0
\put{$\sssize 0$} at 2 0
\endpicture} at 5 2
\put{\beginpicture
\setcoordinatesystem units <.15cm,.1cm>
\put{$\sssize 0$} at 0 1
\put{$\sssize 1$} at 0 -1
\put{$\sssize 1$} at 1 0
\put{$\sssize 0$} at 2 0
\endpicture} at 6 2
\put{\beginpicture
\setcoordinatesystem units <.15cm,.1cm>
\put{$\sssize 0$} at 0 1
\put{$\sssize 0$} at 0 -1
\put{$\sssize 1$} at 1 0
\put{$\sssize 1$} at 2 0
\endpicture} at 7 2
\put{\beginpicture
\setcoordinatesystem units <.15cm,.1cm>
\put{$\sssize 1$} at 0 1
\put{$\sssize 1$} at 0 -1
\put{$\sssize 1$} at 1 0
\put{$\sssize 0$} at 2 0
\endpicture} at 8 2
\put{\beginpicture
\setcoordinatesystem units <.15cm,.1cm>
\put{$\sssize 0$} at 0 1
\put{$\sssize 1$} at 0 -1
\put{$\sssize 1$} at 1 0
\put{$\sssize 1$} at 2 0
\endpicture} at 9 2
\put{\beginpicture
\setcoordinatesystem units <.15cm,.1cm>
\put{$q_B(\bdim N)$} at 0 0
\endpicture} at -1.5 1
\put{\beginpicture
\setcoordinatesystem units <.15cm,.1cm>
\put{$q_{B'}(\bdim N)$} at 0 0
\endpicture} at -1.5 0
\multiput{$1$} at 0 1  1 1  2 1  3 1    5 1  6 1  7 1  8 1  9 1  
   0 0  1 0  2 0  3 0  4 0  5 0   7 0  8 0   /
\multiput{$3$} at 4 1  6 0  9 0 /
\plot -2.5 1.5  9.5 1.5 /
\plot -2.5 0.5  9.5 0.5 /
\plot -0.5 2.3   -0.5 -0.3 /
\endpicture}
$$
	\bigskip
{\bf 13.3.} Let us consider now canonical algebras.
A canonical algebra $C$ is a tilted algebra
if and only if it is domestic (thus, the quiver obtained from the quiver of $C$
by deleting the source is a Dynkin diagram), and these algebras are cluster-concealed.
For example, let us consider the canonical algebra of type $\Bbb E_7$, its quiver has
the form
$$
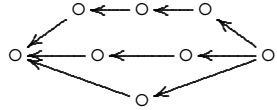
\hbox{\beginpicture
\setcoordinatesystem units <.7cm,.5cm>
\multiput{$\circ$} at 0 0  1 1  2 1  3 1  4 0  1.3 0  2.7 0  2 -1 /
\arr{1.8 1}{1.2 1}
\arr{2.8 1}{2.2 1}

\arr{3.8 0}{2.9 0}
\arr{2.5 0}{1.5 0}
\arr{1.1 0}{0.2 0}

\arr{0.8 0.8}{.2 .2}
\arr{3.8 0.2}{3.2 0.8}

\arr{3.8 -0.1}{2.2 -.9}
\arr{1.8 -.9}{0.2 -0.1}
\endpicture}
$$
and there is a single relation: the sum of the paths from the source to the sink.
The corresponding cluster-tilted algebra has one additional arrow $\gamma$:
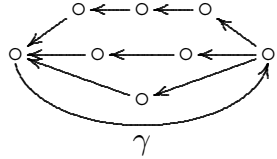
$$
\hbox{\beginpicture
\setcoordinatesystem units <.7cm,.5cm>
\multiput{$\circ$} at 0 0  1 1  2 1  3 1  4 0  1.3 0  2.7 0  2 -1 /
\arr{1.8 1}{1.2 1}
\arr{2.8 1}{2.2 1}

\arr{3.8 0}{2.9 0}
\arr{2.5 0}{1.5 0}
\arr{1.1 0}{0.2 0}

\arr{0.8 0.8}{.2 .2}
\arr{3.8 0.2}{3.2 0.8}

\arr{3.8 -0.1}{2.2 -.9}
\arr{1.8 -.9}{0.2 -0.1}
\ellipticalarc axes ratio 1.8:1 165 degrees from 0 -.2  center at 2 0
\arr{3.97 -0.3}{3.995 -0.2}
\put{$\gamma$} at 2 -2
\endpicture}
$$
and a lot of zero relations: any path from a vertex $x$ to a vertex $y$ and
involving $\gamma$ is a zero relation, provided the quiver of $B$ contains
an arrow $y \to x$.

We consider the indecomposable $B_2$-module $N$ as well as its image
$Z = \pi(N)$ under the covering functor $\pi$:
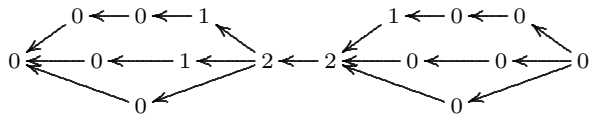
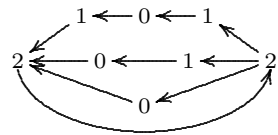
$$
\hbox{\beginpicture
\setcoordinatesystem units <.7cm,.5cm>
\put{\beginpicture

\put{$\bdim N$} at -3 1.5
\put{\beginpicture

\multiput{$\ssize 0$} at  0 0  1 1  2 1    1.3 0    2 -1 /
\multiput{$\ssize 1$} at  3 1   2.7 0   /
\multiput{$\ssize 2$} at  4 0 /

\arr{1.8 1}{1.2 1}
\arr{2.8 1}{2.2 1}

\arr{3.8 0}{2.9 0}
\arr{2.5 0}{1.5 0}
\arr{1.1 0}{0.2 0}

\arr{0.8 0.8}{.2 .2}
\arr{3.8 0.2}{3.2 0.8}

\arr{3.8 -0.1}{2.2 -.9}
\arr{1.8 -.9}{0.2 -0.1}
\endpicture} at 0 0

\arr{2.8 0}{2.2 0}

\put{\beginpicture

\multiput{$\ssize 0$} at  2 1  3 1  1.3 0  2.7 0  2 -1  4 0 /
\multiput{$\ssize 1$} at  1 1     /
\multiput{$\ssize 2$} at 0 0   /

\arr{1.8 1}{1.2 1}
\arr{2.8 1}{2.2 1}

\arr{3.8 0}{2.9 0}
\arr{2.5 0}{1.5 0}
\arr{1.1 0}{0.2 0}

\arr{0.8 0.8}{.2 .2}
\arr{3.8 0.2}{3.2 0.8}

\arr{3.8 -0.1}{2.2 -.9}
\arr{1.8 -.9}{0.2 -0.1}
\endpicture} at 5 0
\endpicture} at 0 0

\put{\beginpicture
\setcoordinatesystem units <.7cm,.5cm>
\put{$\bdim Z$} at -1 1.5

\multiput{$\ssize 0$} at  2 1    1.3 0    2 -1 /
\multiput{$\ssize 1$} at  1 1   3 1   2.7 0   /
\multiput{$\ssize 2$} at 0 0   4 0 /
\arr{1.8 1}{1.2 1}
\arr{2.8 1}{2.2 1}

\arr{3.8 0}{2.9 0}
\arr{2.5 0}{1.5 0}
\arr{1.1 0}{0.2 0}

\arr{0.8 0.8}{.2 .2}
\arr{3.8 0.2}{3.2 0.8}

\arr{3.8 -0.1}{2.2 -.9}
\arr{1.8 -.9}{0.2 -0.1}
\ellipticalarc axes ratio 1.8:1 165 degrees from 0 -.2  center at 2 0
\arr{3.97 -0.3}{3.995 -0.2}
\endpicture} at 10 0

\endpicture}
$$
Note that $q_B(\bdim Z) = 9$. Here, both $B$-modules $N_1,N_2$ are 
decomposable, and it is easy to see that 
$\dim \End(N_1) = 3,$ and $\dim\End(N_2) = 2$;  the $B$-modules
$N_1, N_2$ are the following:

$$
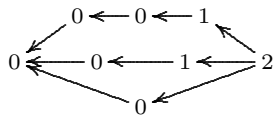
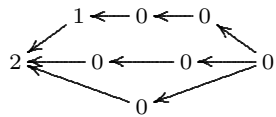
\hbox{\beginpicture
\setcoordinatesystem units <.7cm,.5cm>

\put{\beginpicture

\multiput{$\ssize 0$} at  0 0  1 1  2 1    1.3 0    2 -1 /
\multiput{$\ssize 1$} at  3 1   2.7 0   /
\multiput{$\ssize 2$} at  4 0 /

\arr{1.8 1}{1.2 1}
\arr{2.8 1}{2.2 1}

\arr{3.8 0}{2.9 0}
\arr{2.5 0}{1.5 0}
\arr{1.1 0}{0.2 0}

\arr{0.8 0.8}{.2 .2}
\arr{3.8 0.2}{3.2 0.8}

\arr{3.8 -0.1}{2.2 -.9}
\arr{1.8 -.9}{0.2 -0.1}
\endpicture} at 0 -2

\put{\beginpicture

\multiput{$\ssize 0$} at  2 1  3 1  1.3 0  2.7 0  2 -1  4 0 /
\multiput{$\ssize 1$} at  1 1     /
\multiput{$\ssize 2$} at 0 0   /

\arr{1.8 1}{1.2 1}
\arr{2.8 1}{2.2 1}

\arr{3.8 0}{2.9 0}
\arr{2.5 0}{1.5 0}
\arr{1.1 0}{0.2 0}

\arr{0.8 0.8}{.2 .2}
\arr{3.8 0.2}{3.2 0.8}

\arr{3.8 -0.1}{2.2 -.9}
\arr{1.8 -.9}{0.2 -0.1}
\endpicture} at 5 -2

\put{$\bdim N_1$} at -3 0
\put{$\bdim N_2$} at 3 0

\put{} at 14 0

\endpicture}
$$
	\medskip
We exhibit another indecomposable $B_2$-module $N'$ 
as well as its image
$Z = \pi(N)$ under the covering functor $\pi$:
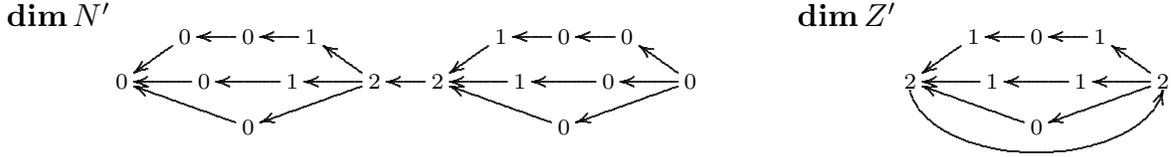
$$
\hbox{\beginpicture
\setcoordinatesystem units <.7cm,.5cm>
\put{\beginpicture

\put{$\bdim N'$} at -3 1.5

\put{\beginpicture

\multiput{$\ssize 0$} at  0 0  1 1  2 1    1.3 0    2 -1 /
\multiput{$\ssize 1$} at  3 1   2.7 0   /
\multiput{$\ssize 2$} at  4 0 /

\arr{1.8 1}{1.2 1}
\arr{2.8 1}{2.2 1}

\arr{3.8 0}{2.9 0}
\arr{2.5 0}{1.5 0}
\arr{1.1 0}{0.2 0}

\arr{0.8 0.8}{.2 .2}
\arr{3.8 0.2}{3.2 0.8}

\arr{3.8 -0.1}{2.2 -.9}
\arr{1.8 -.9}{0.2 -0.1}
\endpicture} at 0 0
\arr{2.8 0}{2.2 0}
\put{\beginpicture

\multiput{$\ssize 0$} at  2 1  3 1  2.7 0  2 -1  4 0 /
\multiput{$\ssize 1$} at  1 1  1.3 0   /
\multiput{$\ssize 2$} at 0 0   /

\arr{1.8 1}{1.2 1}
\arr{2.8 1}{2.2 1}

\arr{3.8 0}{2.9 0}
\arr{2.5 0}{1.5 0}
\arr{1.1 0}{0.2 0}

\arr{0.8 0.8}{.2 .2}
\arr{3.8 0.2}{3.2 0.8}

\arr{3.8 -0.1}{2.2 -.9}
\arr{1.8 -.9}{0.2 -0.1}
\endpicture} at 5 0
\endpicture} at 0 0

\put{\beginpicture
\setcoordinatesystem units <.7cm,.5cm>
\put{$\bdim Z'$} at -1 1.5

\multiput{$\ssize 0$} at  2 1      2 -1 /
\multiput{$\ssize 1$} at  1 1   3 1  1.3 0   2.7 0   /
\multiput{$\ssize 2$} at 0 0   4 0 /
\arr{1.8 1}{1.2 1}
\arr{2.8 1}{2.2 1}

\arr{3.8 0}{2.9 0}
\arr{2.5 0}{1.5 0}
\arr{1.1 0}{0.2 0}

\arr{0.8 0.8}{.2 .2}
\arr{3.8 0.2}{3.2 0.8}

\arr{3.8 -0.1}{2.2 -.9}
\arr{1.8 -.9}{0.2 -0.1}
\ellipticalarc axes ratio 1.8:1 165 degrees from 0 -.2  center at 2 0
\arr{3.97 -0.3}{3.995 -0.2}
\endpicture} at 10 0

\endpicture}
$$
Note that the composition of the horizontal maps on the left,
as well as the composition of the corresponding maps on the right have to be zero.
On the right, we see an indecomposable $B^c$-module $Z'$ such that
there is an arrow (namely the one in the center) 
such that the corresponding vector space map used
in $Z'$ is the zero map and neither an epimorphism nor a monomorphism. 
Here,
$q_B(\bdim Z') = 7.$

	   \bigskip\bigskip
{\bf 14. Tilting modules which are neither preprojective nor preinjective.}
     \medskip
The separation property holds only for preprojective (or preinjective) tilting modules,
as we are going to show now. As above, let $A$ be a finite-dimensional hereditary
$k$-algebra. 
	     \medskip
{\bf Proposition 8.} {\it 
Let $T_i$ be an indecomposable regular $A$-module. Then the component
containing $T_i$ contains infinitely many indecomposable modules $M$ such that both
$\Hom(T_i,M) \neq 0$ and $\Ext^1(T_i,M) \neq 0.$ If $A$ is wild,
then any regular component contains infinitely many indecomposable modules $M$ with this
property.}
	\medskip
Proof: In case $A$ is tame, we deal with a stable tube and the stated property
is easy to check. Thus, assume that $A$ is wild. 
Note that a regular component of a wild hereditary algebra 
is of the form $\Bbb ZA_\infty$,

We use 
the following well-known assertion of Baer [B] and Kerner [K]:
If $X$ and $Y$ are indecomposable regular modules, then
$\Hom_A(X,\tau^mY) \neq 0$ and $\Hom(\tau^{-m} Y,X) \neq 0$
for $m \gg 0$. Thus, consider any regular component of
the Auslander-Reiten quiver of $A$ and let $N$ be an indecomposable 
module in this component which is quasi-simple.
Then,  we have
$\Hom(T_i,\tau^m N) \neq 0$ and $\Hom(\tau^{-m} N,T_i) \neq 0$
for $m\gg 0.$ Take such a natural number $m\ge 1$
and consider the ray starting in $\tau^m N$
$$
 \tau^mN = \tau^mN[1] \to \tau^mN[2] \to \cdots \to \tau^mN[t] \to 
 \cdots,
$$
it consists of indecomposable modules and irreducible monomorphisms. Let 
$$
 M = \tau^mN[2m].
$$
Now $M$ has a filtration going upwards with
factors $\tau^mN,\tau^{m-1}N,\dots, \tau^{-m+1}N$. In particular, it follows that 
$M$ maps onto $\tau^{-m+1}N$.
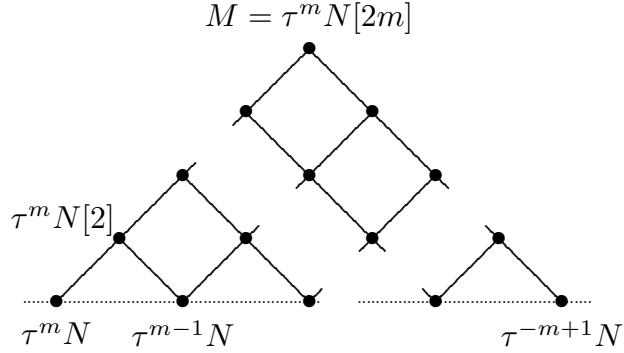
$$
\hbox{\beginpicture
\setcoordinatesystem units <.7cm,.7cm>
\multiput{$\bullet$} at 0 0  1 1  2 2  2 0  3 1  4 0  3 3  4 4  4 2  5 3  5 1  6 2
     6 0  7 1  8 0  /
\multiput{} at 0 0  8 5 /
\plot 0 0  2.2 2.2 /
\plot 1 1  2 0  3.2 1.2 /
\plot 2 2  4 0 /
\plot 3 3  5.2 0.8 /
\plot 2.8 2.8  4 4  6.2 1.8 /
\plot 3.8 1.8  5 3 /
\plot 4.8 0.8  6 2 /
\plot 4 0  4.2 0.2 /
\plot 5.8 0.2  6 0 /
\plot 6 0  7 1  /
\plot 6.8 1.2  8 0 /
\setdots <0.5mm>
\plot  -.5 0  4.2 0 /
\plot 4.8 0  8.5 0 / 
\put{$\tau^mN$} at 0 -0.5
\put{$\tau^{m-1}N$} at 2 -0.5
\put{$M = \tau^mN[2m]$} at 4 4.5
\put{$\tau^mN[2]$} at 0.1 1.3
\put{$\tau^{-m+1}N$} at 8 -0.5

\endpicture}
$$

We claim that both $\Hom(T_i,M)$   and $\Ext^1(T_i,M)$ are non-zero.
On the one hand, $0\neq \Hom(T_i,\tau^mN)$ embeds into $\Hom(T_i,M)$,
since $\tau^mN$ is a submodule of $M$. This shows that $\Hom(T_i,M)\neq 0.$ 
On the other hand,  
we have $0 \neq \Hom(\tau^{-m}N,T_i) \simeq \Hom(\tau^{-m+1}N,\tau T_i),$
since $\tau^{-1}$ is left adjoint to $\tau$.
Composing a surjective map 
$M \to  \tau^{-m+1}N$ with 
a non-zero map $\tau^{-m+1}N \to \tau T_i$, we obtain a non-zero map. This shows that
$\Hom(M,\tau T_i) \neq 0,$ and therefore 
$$
 \Ext^1(T_i,M) \simeq D\Hom(M,\tau T_i) \neq 0.
$$
This completes the proof.

	\bigskip
{\bf Proposition 9.} {\it 
If $T_0$ is an indecomposable preprojective module, 
and $T_\infty$ is an indecomposable preinjective module, 
then
any regular component contains infinitely many indecomposable modules $M$
with $\Hom(T_0,M) \neq 0$ and $\Ext^1(T_\infty,M) \neq 0.$}
     \medskip
Proof: The proof is similar to the previous proof. Here, we use that for indecomposable
modules $P,R,Q$ with $P$ preprojective, $R$ regular and $Q$ preinjective,
$\Hom(P,\tau^m R) \neq 0$ and $\Hom(\tau^{-m} R,Q) \neq 0$ for $m \gg 0$.
	       \bigskip
{\bf Corollary.} {\it
 Let $A$ be a finite-dimensional hereditary algebra and
$T$ a tilting module. 
The following conditions are equivalent:
\item{\rm(i)} The tilting module $T$ is neither preprojective nor preinjective.
\item{\rm(ii)} $\Cal M(T)$ is infinite.
\item{\rm(iii)} $\Cal M(T)$ contains a regular module.
\item{\rm(iv)} $\Cal M(T)$ contains indecomposable 
  regular modules of arbitrarily large length.}
  \medskip
Proof: Of course, (iv) implies both (ii) and (iii).

Now assume that (ii) or (iii) holds. 
If $T$ is preprojective, then we have seen in Lemma 3
that $\Cal M(T)$ is a finite set of preprojective modules; similarly, if $T$
is preinjective, then $\Cal M(T)$ is a finite set of preinjective modules.
This contradiction shows that $T$ cannot be preprojective or preinjective, thus (i)
holds. 

Conversely, let us assume (i).
Then either $T$ has an indecomposable summand which is regular, or two 
indecomposable summands one being preprojective, the other being preinjective.
If $A$ is wild, then
the previous two propositions show that $\Cal M(T)$ contains infinitely
many indecomposable modules which belong to the same regular component. But
a regular component cf a hereditary algebra $A$ 
contains only finitely many indecomposable modules of a given
length (for $A$ wild, see for example [Zg]). 
Actually, one easily observes that the proofs
of the two previous propositions yield sequences of indecomposable modules in
$\Cal M(T)$ of unbounded length.) This shows (iv) in case $A$ is wild.
In case $A$ is tame, the proof is similar. 
   \bigskip
{\bf 15. A further remark.}
     \medskip
Starting with a finite-dimensional hereditary algebra $A$ and a tilting $A$-module
$T$, one considers $T$ as an object in $\Cal C_A$ and obtains in this way 
a so-called cluster-tilting object. 
However, one knows that not all cluster-tilted objects
of $\Cal C_A$ are obtained in this way --- it may be necessary to change the
orientation of the quiver of $A$. For the benefit of the reader let us include
an easy recipe for finding an orientation such that a given cluster-tilting object
can be considered as a module. Of course, we only have to consider the cases
when $T$ is not regular.
     \medskip
{\bf Proposition 10.} {\it Let $T$ be a cluster-tilting object in a cluster
category $\Cal C$.
Let $\Cal S$ be a slice in $\Cal C$ such that 
the sources of $\Cal S$ belong to $\add T$. Then no indecomposable direct summand of $T$
belongs to $\tau^{-1} \Cal S$.}
	\medskip
Proof: Assume $T'$ is an indecomposable direct summand of $T$ and 
belongs to $\tau^{-1} \Cal S$. Then $\tau T' \in \Cal S$. There is
a source $S$ in $\Cal S$ with $\Hom(S,\tau T')\neq 0$.
Thus $\Ext^1(T',S) = D\Hom(S,\tau T') \neq 0,$ thus $S$ cannot be a direct
summand of $T$, in contrast to the assumption.
	\medskip
We can apply the proposition as follows: Let $T_1$ be any indecomposable
direct summand of $T$ which is not regular. Then there is a unique slice 
in $\Cal C$ such that $T_1$ is the unique source.

   \bigskip\bigskip
{\bf References.}
     \medskip
\item{[ABS]} Assem, Br\"ustle and Schiffler: Cluster-tilted algebras as trivial
  extensions. To appear in J. London Math. Soc.
 
\item{[B]} Baer, D.: Wild hereditary artin algebras and linear methods.
  manuscripta math. 55 (1986), 69-82.
\item{[BMR]} Buan, Marsh, Reiten: Cluster tilted algebras. 
  Trans. Amer. Math. Soc., 359, no. 1, 323--332 (2007) 
\item{[BRS]} Buan, Reiten, Seven: Tame concealed algebras and cluster
  quivers of minimal infinite type. J. Pure Applied Algebra 211 (2007), 71-82.
\item{[D]} Drozd: Matrix problems and categories of matrices. Zap. Nauchn. Sem.
   Leningrad. Otdel Mat. Inst. Steklov (LOMI) 28 (1972), 144-153.
\item{[DLS]} Dowbor, Lenzing, Skowronski: Galois coverings of algebras by
 locally support-finite categories. In: Representation Theory I. Finite Dimensional
  Algebras (ed. Dlab, Gabriel, Michler). Springer LNM 1177 (1986), 91-93.
\item{[GP]} Geng, Peng: The dimension vectors of indecomposable modules of
  cluster-tilted algebras and the Fomin-Zelevinsky
  denominators conjecture. Preprint (2009).
\item{[HV]} Happel-Vossieck: Minimal algebras of infinite representation type with 
  preprojective component. manuscripta math. 42 (1983), 221-243.
\item{[KR]} Keller-Reiten: Cluster-tilted algebras are Gorenstein
   and stably Calabi-Yau. Adv. Math. 211 (2007), 123-151
\item{[K]} Kerner, O.: Tilting wild algebras. J. London Math. Soc. 39 (1989), 29-47. 
\item{[P]} Parsons, M.: On indecomposable modules over cluster-tilted algebras
of type $A$. 
PhD Thesis, University of Leicester, 2007.
\item{[PS]} de la Pe\~na, Simson: Prinjective modules, reflection functors,
 quadratic forms and Auslander-Reiten sequences. Trans. Amer. Math. Soc. 329 (1992), 
  733-753.
\item{[R]} Ringel: Some remarks concerning tilting modules and
 tilted algebra. Origin. Relevance. Future. Appendix to the Handbook of Tilting
 Theory (ed. Angeleri-H\"ugel, Happel, Krause). London Math. Soc. LNS 332 (2007), 413-472.
\item{[Si]} Simson: Linear representations of partially ordered sets
    and vector spaces categories. Gordon and Breach (1992).
\item{[Se]} Seven: Recognizing cluster algebras of finite type, 
  Electronic J Combinatorics. 14 (2007), 35 pp.
\item{[U]} Unger: The concealed algebras of the minimal wild, hereditary algebras
  Bayreuth. Math. Schr, 1990.
\item{[Z]} Zhu, Bin: Equivalences between cluster categories. 
  Journal of Algebra, Vol.304, 832-850,2006 
\item{[Zg]} Zhang, Y.: The modules in any component of the AR-quiver of a 
 wild hereditary artin algebra are uniquely determined by their composition 
 factors. Archiv Math. 53 (1989), 250-251. 
 \bigskip\bigskip
 \bigskip\bigskip
{\bf Appendix: Torsion pairs and matrix categories.}
     \medskip
For the convenience of the reader, we add a proof of the  result
which seems to be folklore. Assertions of this kind can be traced back to
the Kiev school of Nazarova, Roiter and Drozd.
    \medskip
{\bf Proposition.} {\it 
Let $(\Cal F,\Cal G)$ be a torsion pair in the abelian category $\Cal A$.
Given $A\in \Cal A$, let $\epsilon_A$
be the equivalence class of the canonical exact sequence 
$0 \to tA \to A \to A/tA \to 0$, and $\eta(A) = (A/tA,tA,\epsilon_A)$. Then this
defines a functor $\eta\: \Cal A \to \Mat\Ext^1(\Cal F,\Cal G)$ which 
is full and dense and its kernel is the ideal generated by all 
maps $\Cal F \to \Cal G$.}
     \medskip
In particular, in case $\Cal A$ is a Krull-Remak-Schmidt category, then we see
that the kernel of the functor $\eta$ lies in the radical of $\Cal A$ and therefore
$\eta$ is a representation equivalence. 
       \medskip
Proof: Denote the inclusion maps $tA \to A$ and $tA' \to A'$ by $u,u'$,
respectively, and the projection maps $A \to A/tA$ and $A' \to A'/tA'$
by $p,p'$ respectively.
Let $\alpha\:A \to A'$ be a morphism in $\Cal A.$ Then $\alpha$ maps
$tA$ into $tA'$, thus it induces a map $t\alpha\:tA \to tA'$ as well as
a map $\overline\alpha\:A/tA \to A'/tA',$  thus $u \alpha = (t\alpha) u',$
and $p\overline\alpha  = \alpha p'.$

Using $t\alpha$, we obtain the following induced exact sequence
$$
\CD
 0 @>>> tA @>u>>      A @>p>> A/tA @>>> 0 \cr
 @.   @Vt\alpha VV  @V\alpha' VV   @|    \cr
 0 @>>> tA' @>>v >      B @>>q> A/tA @>>> 0 .
\endCD
$$
Since $u\alpha = (t\alpha)u'$, there is a 
map $\alpha''\:B \to A'$ such that $\alpha = \alpha'\alpha''$
and $v\alpha'' = u'$. It follows that $q\overline\alpha  = \alpha''p'$,
since 
$$
 \alpha'\alpha''p' = \alpha p' = p\overline\alpha = \alpha'q\overline\alpha,
$$
and 
$$
 v\alpha''p' = u'p' = 0 = vq\overline\alpha.
$$

Thus we also have a commutative diagram with exact rows:
$$
\CD
 0 @>>> tA' @>v >>      B @>q>> A/tA @>>> 0  \cr
 @.   @|  @V\alpha'' VV   @V\overline\alpha VV    \cr
 0 @>>> tA' @>>u'>      A' @>>p'> A'/tA' @>>> 0.
\endCD
$$
The diagram shows that the upper row is induced from the lower one 
by $\overline \alpha,$ therefore 
$t\alpha \epsilon_A  = \epsilon_{A'}\overline \alpha .$
Thus, we see that $\eta$ is a functor.
      \medskip
First, we show that $\eta$ is dense: any element $\epsilon\in \Ext^1(F,G)$
is the equivalence class of an exact sequence
$$
 0 @>>> G @>\mu >> A @>\epsilon>> F @>>> 0.
$$
Here, the image of $\mu$ has to be $tA$, let $u\:tA \to A$
be the inclusion map and $A/tA$ the canonical projection.
We obtain a commutative diagram with exact rows:
$$
\CD
 0 @>>> G @>\mu>>  A @>\epsilon>> F @>>> 0 \cr
 @.   @V\mu'VV         @|   @VV\epsilon'V              \cr
 0 @>>> tA @>>v >  A @>>q>      A/tA @>>> 0 
\endCD
$$
where $\mu'$ and $\epsilon'$ are isomorphisms. 

      \medskip
Next, we show that $\eta$ is full. Let $A,A'$ be objects and assume
that there are given maps $\beta\:tA \to tA'$ and
$\gamma\:A/tA \to A'/tA'$ such that 
 $\beta \epsilon_A  = \epsilon_{A'}\gamma .$ We obtain the following
commutative diagram with exact rows:
$$
\CD
 0 @>>> tA @>u>>      A @>p>> A/tA @>>> 0 \cr
 @.   @V\beta VV  @V\beta' VV   @|    \cr
 0 @>>> tA' @>>v >      B @>>q> A/tA @>>> 0 \cr
 @.   @|  @V\delta VV   @|    \cr
 0 @>>> tA' @>v' >>      B' @>q'>> A/tA @>>> 0  \cr
 @.   @|  @V\gamma' VV   @V\gamma VV    \cr
 0 @>>> tA' @>>u'>      A' @>>p'> A'/tA' @>>> 0,
\endCD
$$
The upper part shows that the second row is induced from the
first. the lower part shows that the third row is induced from
the forth. The central part means that the two induced sequences are
equivalent: Altogether we obtain the map $\alpha = \beta'\delta\gamma'\:
A \to A'$ and we have $t\alpha = \beta,$ and $\overline\alpha = \gamma$.
Thus, $(\beta,\gamma) = \eta(\alpha).$
      \medskip
It is clear that the maps $\Cal F \to \Cal G$ are in the kernel
of $\eta$. Conversely, assume that $\alpha\:A \to A'$
is in the kernel of $\eta$, thus $t\alpha = 0$ and $\overline \alpha = 0.$
Then $\alpha = p\alpha'u'$ for some $\alpha'\:A/tA \to tA'$,
thus $\alpha$ lies in the ideal generated by the maps $\Cal F \to
\Cal G.$
     \bigskip\bigskip

{\rmk Fakult\"at f\"ur Mathematik, Universit\"at Bielefeld \par
POBox 100\,131, \ D-33\,501 Bielefeld, Germany \par
e-mail: \ttk ringel\@math.uni-bielefeld.de \par}

\bye